\documentclass[10pt,letterpaper]{article}
\usepackage{blindtext,titlefoot}
\usepackage{authblk}
\usepackage{marvosym}
\usepackage[utf8]{inputenc}
\usepackage[english]{babel}
\usepackage{amsmath}
\usepackage{amsfonts}
\usepackage{amssymb}
\usepackage{makeidx}
\usepackage{graphicx}
\usepackage{lmodern}
\usepackage{kpfonts}
\usepackage{dsfont}
\usepackage{cancel}
\usepackage{amsthm} 
\usepackage[left=2cm,right=2cm,top=2cm,bottom=2cm]{geometry}
\usepackage{subcaption}
\usepackage{multirow}
\usepackage{float}
\usepackage{url}
\usepackage{breakurl}
\usepackage[breaklinks]{hyperref}
\usepackage{multicol}
\usepackage{mathtools, cuted}
\usepackage{lipsum}
\usepackage{booktabs}
\usepackage{comment}
\usepackage{cite}

\usepackage{fancyhdr}
\fancypagestyle{firstpage}{\fancyhf{}
\setcounter{page}{1}
\rfoot{\thepage}
}

\providecommand{\nset}[1]{
\mathbb{#1}
}
\providecommand{\set}[1]{
\left\{#1\right\}
}

\providecommand{\ifr}[5]{
{}^{#1}_{#2}{#3}_{#4}^{#5}
}
\providecommand{\gam}[1]{
\Gamma\left(#1 \right)
}

\providecommand{\norm}[1]{
\left\lVert #1 \right\rVert
}
\providecommand{\abs}[1]{
\left\lvert #1 \right\rvert
}
\providecommand{\ds}[1]{
\displaystyle #1
}
\providecommand{\der}[3]{
\dfrac{#1^{#3} }{ #1 #2^{#3}}
}

\DeclareMathOperator{\cond}{cond}
\newtheorem{theorem}{ Theorem}[section]
\newtheorem{definition}[theorem]{Definition}
\newtheorem{proposition}[theorem]{Proposition}

\newtheorem{example}[theorem]{Example}

\usepackage{enumitem} 
\setlist[itemize]{noitemsep} 

\usepackage{titlesec} 
\titleformat{\section}[block]{\large\bfseries\scshape\centering}{\thesection.}{1em}{} 
\titleformat{\subsection}[block]{\large\bfseries\scshape\centering}{\thesubsection.}{1em}{}
\titleformat{\subsubsection}[block]{\large\bfseries\scshape\centering}{\thesubsubsection.}{1em}{} 

\title{\bfseries \huge  Numerical solution using radial basis functions for multidimensional fractional partial differential equations of type Black-Scholes}

\author[,a]{A. Torres-Hernandez  \footnote{Email address: anthony.torres@ciencias.unam.mx}}
\affil[a]{\normalsize Department of Physics, Faculty of Science - UNAM, Mexico}

\author[,b]{F. Brambila-Paz \footnote{Email address: fernandobrambila@gmail.com}}
\affil[b]{\normalsize Department of Mathematics, Faculty of Science - UNAM, Mexico}

\author[,c]{C. A. Torres-Martínez \footnote{Email address: inocencio3@gmail.com}}
\affil[c]{\normalsize Department of Mathematics - UACM, Mexico}

\date{}

\begin{document}

\maketitle


\begin{abstract}

The aim of this paper is to solve numerically, using the meshless method via radial basis functions,  time-space-fractional partial differential equations of type Black-Scholes. The time-fractional partial differential equation appears in several diffusion problems used in physics and engineering applications, and models subdiffusive and superdiffusive behavior of the prices at the stock market. This work shows the flexibility of the radial basis function scheme to solve multidimensional problems with several types of nodes and it also shows how to reduce the condition number of the matrices involved.

\textbf{Keywords}: Fractional differential equations, Meshless methods, Radial Basis Functions, Black-Scholes equations.
\end{abstract}

\section{Introduction}

A fractional derivative is an operator that generalizes the ordinary derivative, in the sense that if

\begin{eqnarray*}
\dfrac{d^\alpha}{d x^\alpha},
\end{eqnarray*}

denotes the differential of order $ \alpha $, it can take values $\alpha\in \nset{R}$  and the first derivative is the particular case when $\alpha=1$. On the other hand, a fractional differential equation is an equation that involves at least one differential operator of order $ \alpha $ with $(n-1)\leq \alpha < n$, for some positive integer $n$, and it is said to be a differential equation of order $\alpha$ if this operator is the highest order in the equation.

The growing interest in fractional calculus has been motivated by applications of fractional equations in different areas of research such as magnetic field theory, fluid dynamics, electrodynamics, multidimensional processes, etc. One of the most popular examples is the convection-diffusion equations, the solution of which can be interpreted as a probability distribution of one or more underlying stochastic processes. One of the most popular examples is the convection-diffusion equations \cite{barkai2000continuous,blumen1989transport,chaves1998fractional,piryatinska2005models}, in which the solutions may be interpreted as a probability distribution of one or more underlying stochastic processes \cite{safdari2015radial}.

The applications of fractional operators have spread to other fields such as finance \cite{safdari2015radial, sabatelli2002waiting} and also in the study for the manufacture of hybrid solar receivers \cite{torres2020reduction}. It should be mentioned that there is also a growing interest in fractional operators and their properties for the solution of nonlinear systems \cite{gdawiec2020newton,torreshern2020,cordero2019variant}. Stochastic processes in financial mathematics may be modeled using Wiener processes or Brownian motion, leading to diffusion partial differential equations. But, if the stochastic process is heavy-tailed rather than Gaussian, then the governing equations are fractional partial differential equations \cite{scalas2000fractional}.

Considering the following partial differential equation, which corresponds to a Black-Scholes model (with a source term $f_I$) and whose details, as well as its deduction, can be found in the reference \cite{goetzmann2014modern}

\begin{eqnarray}\label{eq:S1-01}
\left\{
\begin{array}{cc}
\dfrac{\partial}{\partial \tau}f(S,\tau)+\dfrac{1}{2}\widetilde{\sigma}^2 S^2 \dfrac{\partial^2}{\partial S^2}f(S,\tau)+\widetilde{r}S\dfrac{\partial}{\partial S}f(S,\tau)-\widetilde{r}f(S,\tau)=f_I(S,\tau), &(S,\tau)\in \widetilde{\Omega} \times \widetilde{ D}\\
f(S,\tau)=f_B(S,\tau), & (S,\tau)\in \partial \widetilde{\Omega} \times \widetilde{D}\\
f(S,\tau_0)=f_0(S), & S\in \widetilde{\Omega}
\end{array}\right.,
\end{eqnarray}

with $\widetilde{\Omega}$ and $\widetilde{D}$ subsets of $\nset{R}_{\geq 0}$. It should be mentioned that a complete study of the Black-Sholes model goes beyond the purpose of this document, our interest will focus only on finding the numerical solution of some variations of the previous model with fractional operators. Considering $\widetilde{D}$ a finite interval and  using the change of variables

\begin{align*}
\tau=&t_m-t,\\
S=&e^x,
\end{align*}

we obtain that

\begin{align*}
\dfrac{\partial}{\partial \tau}=&\left( \dfrac{\partial \tau}{\partial t}\right)^{-1}\dfrac{\partial}{\partial t}=-\dfrac{\partial}{\partial t},\\
\dfrac{\partial}{\partial S}=&\left( \dfrac{\partial S}{\partial x}\right)^{-1}\dfrac{\partial}{\partial x}=e^{-x}\dfrac{\partial}{\partial x},\\
\dfrac{\partial^2}{\partial S^2}=&e^{-x}\dfrac{\partial}{\partial x}\left(e^{-x}\dfrac{\partial}{\partial x} \right)=e^{-2x}\dfrac{\partial^2}{\partial x^2}-e^{-2x}\dfrac{\partial}{\partial x},
\end{align*}

as a consequence

\begin{align*}
\dfrac{\partial}{\partial \tau}f(S,\tau)=&-\dfrac{\partial}{\partial t}f(e^x,t_m-t)=-\dfrac{\partial}{\partial t}u(x,t),\\
S\dfrac{\partial}{\partial S}f(S,\tau)=&e^{x}\left( e^{-x}\dfrac{\partial}{\partial x}\right)f(e^x,t_m-t)=\dfrac{\partial}{\partial x}u(x,t),\\
S^2\dfrac{\partial^2}{\partial S^2}f(S,t)=&e^{2x}\left(e^{-2x}\dfrac{\partial^2}{\partial x^2}-e^{-2x}\dfrac{\partial}{\partial x}\right)f(e^x,t_m-t)=\dfrac{\partial ^2}{\partial x^2}u(x,t)-\dfrac{\partial}{\partial x}u(x,t),
\end{align*}

therefore it is possible to rewrite \eqref{eq:S1-01} as follows

\begin{eqnarray}\label{eq:S1-02}
\left\{
\begin{array}{cc}
\dfrac{\partial}{\partial t}u(x,t)-\dfrac{1}{2}\widetilde{\sigma}^2  \dfrac{\partial^2}{\partial x^2}u(x,t)-\left(\widetilde{r}-\dfrac{1}{2}\widetilde{\sigma}^2  \right) \dfrac{\partial}{\partial x}u(x,t)+\widetilde{r}u(x,t)=u_I(x,t), &(x,t)\in \Omega \times D\\
u(x,t)=u_B(x,t), & (x,t)\in \partial \Omega \times D\\
u(x,t_0)=u_0(x), & x\in \Omega
\end{array}\right. .
\end{eqnarray}

The above equation may be generalized considering fractional operators and larger dimensions using the following expression

\begin{eqnarray}\label{eq:S1-03}
\left\{
\begin{array}{cc}
\dfrac{\partial^\alpha}{\partial t^\alpha}u(x,t)-\mathcal{L}_{\beta,r}u(x,t)=u_I(x,t), & (x,t)\in \Omega \times D\\
u(x,t)=u_B(x,t), & (x,t)\in \partial \Omega \times D\\
u(x,t_0)=u_0(x), & x\in \Omega
\end{array}\right. ,
\end{eqnarray}

with

\begin{eqnarray}\label{eq:S1-04}
\mathcal{L}_{\beta,r}:= \dfrac{1}{2}\widetilde{\sigma}^2  \dfrac{\partial^{\beta+1}}{\partial r^{\beta+1}}+\left(\widetilde{r}-\dfrac{1}{2}\widetilde{\sigma}^2  \right) \dfrac{\partial^\beta}{\partial r^\beta}-\widetilde{r},
\end{eqnarray}

where $0<\alpha,\beta\leq 1$ and $r=\norm{x}_2$ with $x\in \nset{R}^d$. It should be noted that when $\alpha=\beta=d=1$, the equation \eqref{eq:S1-03} coincides with the equation \eqref{eq:S1-02}. In the following sections, the parts necessary to find the numerical solution of the equation \eqref{eq:S1-03} will be given in as much detail as possible.

\section{Basic Definitions of the Fractional Derivative}

\subsection{Introduction to the Definition of Riemann-Liouville }

One of the key pieces in the study of fractional calculus is the iterated integral, which is defined as follows \cite{hilfer00}

\begin{definition}
Let $ L_{loc} ^ 1 (a, b) $, the space of locally integrable functions in the interval $ (a, b) $. If $ f $ is a function such that $ f \in L_ {loc} ^ 1 (a, \infty) $, then the $n$-th iterated integral of the function $ f $ is given by 

\begin{eqnarray}\label{eq:c1.16}
\begin{array}{c}
\ds \ifr{}{a}{I}{x}{n} f(x)=\ifr{}{a}{I}{x}{}\left(\ifr{}{a}{I}{x}{n-1} f(x)  \right)=\frac{1}{(n-1)!}\int_a^x(x-t)^{n-1}f(t)dt,
\end{array}
\end{eqnarray}

where

\begin{eqnarray*}
\ifr{}{a}{I}{x}{} f(x):=\int_a^x f(t)dt.
\end{eqnarray*}

\end{definition}

Considerate that $ (n-1)! = \gam{n} $
, a generalization of \eqref{eq:c1.16} may be obtained for an arbitrary order $ \alpha> 0 $

\begin{eqnarray}\label{eq:c1.17}
\ifr{}{a}{I}{x}{\alpha} f(x)=\dfrac{1}{\gam{\alpha}}\int_a^x(x-t)^{\alpha-1}f(t)dt,
\end{eqnarray}

similarly, if $ f \in L_{loc} ^ 1 (- \infty, b) $, we may define

\begin{eqnarray}\label{eq:c1.18}
\ifr{}{x}{I}{b}{\alpha} f(x)=\dfrac{1}{\gam{\alpha}}\int_x^b(t-x)^{\alpha -1}f(t)dt,
\end{eqnarray} 

the equations \eqref{eq:c1.17} and \eqref{eq:c1.18} correspond to the definitions of \textbf{right and left fractional integral of Riemann-Liouville}, respectively. The fractional integrals fulfill the  \textbf{semigroup property}, which is given in the following proposition \cite{hilfer00}

\begin{proposition}
Let $ f $ be a function. If $ f \in L_{loc} ^ 1 (a, \infty) $, then the fractional integrals of $ f $ fulfill that

\begin{eqnarray}\label{eq:c1.19}
\ifr{}{a}{I}{x}{\alpha} \ifr{}{a}{I}{x}{\beta}f(x) = \ifr{}{a}{I}{x}{\alpha + \beta}f(x),& \alpha,\beta>0.
\end{eqnarray}

\end{proposition}

From the previous result, and considering that the operator $ d / dx $  is the inverse operator to the left of the operator $ \ifr {}{a}{I}{x}{} $, any integral $ \alpha$-th of a function $ f \in L_{loc} ^ 1 (a, \infty) $ may be written as

\begin{eqnarray}\label{eq:c1.20}
\ifr{}{a}{I}{x}{\alpha}f(x)=\dfrac{d^n}{dx^n}\left( \ifr{}{a}{I}{x}{n}\ifr{}{a}{I}{x}{\alpha}f(x) \right)=\dfrac{d^n}{dx^n}\left( \ifr{}{a}{I}{x}{n+\alpha}f(x)\right).
\end{eqnarray}
 
Considering \eqref{eq:c1.17} and \eqref{eq:c1.20}, we can built the operator  \textbf{Fractional Derivative of Riemann-Liouville} $\ifr{}{a}{D}{x}{\alpha}$, as follows \cite{hilfer00,kilbas2006theory}

\begin{eqnarray}\label{eq:c1.23}
\normalsize
\begin{array}{c}
\ifr{}{a}{D}{x}{\alpha}f(x) := \left\{
\begin{array}{cc}
\ds \ifr{}{a}{I}{x}{-\alpha}f(x), &\mbox{if }\alpha<0\\  
\ds \dfrac{d^n}{dx^n}\left( \ifr{}{a}{I}{x}{n-\alpha}f(x)\right), & \mbox{if }\alpha\geq 0
\end{array}
\right. ,
\end{array}
\end{eqnarray}

where  $ n = \lceil \alpha \rceil$,  then applying the  operator \eqref{eq:c1.23} to the  function $ x^{\mu} $, with  $\alpha \in \nset{R}\setminus\nset {Z} $ and $\mu>-1$, we obtain the following result

\begin{eqnarray}\label{eq:c1.13}
\ifr{}{0}{D}{x}{\alpha}x^\mu = 
 \dfrac{\gam{\mu+1}}{\gam{\mu-\alpha+1}}x^{\mu-\alpha}.
\end{eqnarray}

\subsection{Introduction to the Definition of Caputo }

Michele Caputo (1969) published a book and introduced a new definition of fractional derivative, he created this definition with the objective of modeling anomalous diffusion phenomena. The definition of Caputo had already been discovered independently by Gerasimov (1948). This fractional derivative is of the utmost importance since it allows us to give a physical interpretation of the initial value problems, moreover to being used to model fractional time. In some texts, it is known as the fractional derivative of Gerasimov-Caputo.

Let $ f $ be a function, such that $ f $ is $ n$-times differentiable with $ f ^{(n)} \in L_{loc}^ 1 (a, b) $, then the \textbf{(right) fractional derivative  of Caputo} is defined as \cite{kilbas2006theory}

\begin{align}
\ifr{C}{a}{D}{x}{\alpha}f(x):= &\ifr{}{a}{I}{x}{n-\alpha}\left( \der{d}{x}{n} f(x)\right) = \dfrac{1}{\gam{n-\alpha}}\int_{a}^{x} (x-t)^{n-\alpha -1} f^{(n)}(t)dt , \label{eq:c1.25}
\end{align}

where $n=\lceil \alpha \rceil$. It should be mentioned that the fractional derivative of Caputo behaves as the inverse operator to the left of fractional integral of Riemann-Liouville , that is,

\begin{eqnarray*}
\ifr{C}{a}{D}{x}{\alpha}(\ifr{}{a}{I}{x}{\alpha}f(x))=f(x).
\end{eqnarray*}

On the other hand, the relation between the fractional derivatives of Caputo and Riemann-Liouville is given by the following expression \cite{kilbas2006theory}

\begin{eqnarray*}
\ifr{C}{a}{D}{x}{\alpha}f(x)=\ifr{}{a}{D}{x}{\alpha}\left(f(x)-\sum_{k=0}^{n-1}\dfrac{f^{(k)}(a)}{k!}(x-a)^k\right), 
\end{eqnarray*}

then, if $f^{(k)}(a)=0 \ \ \forall k<n$, we obtain

\begin{eqnarray*}
\ifr{C}{a}{D}{x}{\alpha}f(x)=\ifr{}{a}{D}{x}{\alpha}f(x), 
\end{eqnarray*}

considering the previous particular case, it is possible to unify the definitions of fractional integral of Riemann-Liouville and  fractional derivative of Caputo as follows

\begin{eqnarray}\label{eq:c1.233}
\begin{array}{c}
\ifr{C}{a}{D}{x}{\alpha}f(x) := \left\{
\begin{array}{cc}
\ds \ifr{}{a}{I}{x}{-\alpha}f(x), &\mbox{if }\alpha<0\\  
\ds \ifr{}{a}{I}{x}{n-\alpha}\left( \der{d}{x}{n} f(x)\right) , & \mbox{if }\alpha\geq 0
\end{array}
\right. .
\end{array}
\end{eqnarray}

\subsection{Discretization of the Fractional Derivative of Caputo}

We begin this subsection by considering a uniform partition of the interval $[a,t]$, that is,

\begin{eqnarray*}
a=t_0<t_1<\cdots <t_{m-1}<t_{m}=t,
\end{eqnarray*}

with

\begin{eqnarray*}
t_k=t_0+kdt, & \forall k\geq 0,
\end{eqnarray*}

then, the fractional derivative of Caputo with $(n-1)<\alpha\leq n$ may be written as

\begin{eqnarray*}
\ifr{C}{a}{D}{t}{\alpha}f(t)= \dfrac{1}{\gam{n-\alpha}}\int_{a}^{t} (t-x)^{n-\alpha-1 } f^{(n)}(x)dx
  = \dfrac{1}{\gam{n-\alpha}}\sum_{k=0}^{m-1} \int_{t_{m-k-1}}^{t_{m-k}} (t_{m}-x)^{n-\alpha-1 } f^{(n)}(x)dx,
\end{eqnarray*}

as a consequence

\begin{align}\label{eq:S2-031}
\ifr{C}{a}{D}{t}{\alpha}f(t)    =& \dfrac{1}{\gam{n-\alpha}}\sum_{k=0}^{m-1} \int_{t_{m-k-1}}^{t_{m-k}} (t_{m}-x)^{n-\alpha-1 } \left[ \dfrac{f^{(n-1)}(t_{m-k})-f^{(n-1)}(t_{m-k-1})}{ t_{m-k}-t_{m-k-1}  }+\mathcal{O}( t_{m-k}-t_{m-k-1})   \right]   dx \nonumber \\
  =& \dfrac{1}{\gam{n-\alpha}}\sum_{k=0}^{m-1}\left[\dfrac{(t_{m}-x)^{n-\alpha}}{n-\alpha} \right]_{t_{m-k}}^{t_{m-k-1}} \left[ \dfrac{f^{(n-1)}(t_{m-k})-f^{(n-1)}(t_{m-k-1})}{dt}+\mathcal{O}( dt)   \right] \nonumber \\
    =& \dfrac{1}{\gam{n-\alpha}}\sum_{k=0}^{m-1}\left[\dfrac{(k+1)^{n-\alpha}-k^{n-\alpha}}{n-\alpha}dt^{n-\alpha} \right] \left[ \dfrac{f^{(n-1)}(t_{m-k})-f^{(n-1)}(t_{m-k-1})}{dt}+\mathcal{O}(dt)   \right] \nonumber \\
=& \dfrac{dt^{n-\alpha-1}}{\gam{n-\alpha+1}}\sum_{k=0}^{m-1} \left[(k+1)^{n-\alpha}-k^{n-\alpha} \right]\left[ f^{(n-1)}( t_{m-k}) -f^{(n-1)}( t_{m-k-1} ) \right]+\mathcal{O}\left(
 dt^{n-\alpha+1} \right), 
\end{align}

considering the notation

\begin{eqnarray}\label{eq:S2-032}
c_{\alpha,k}:=(k+1)^{n-\alpha}-k^{n-\alpha}, & n=\lceil \alpha \rceil,
\end{eqnarray}

the equation \eqref{eq:S2-031} may be rewritten as

\begin{align}\label{eq:S2-03}
\ifr{C}{a}{D}{t}{\alpha}f(t)    =&  \dfrac{dt^{n-\alpha-1}}{\gam{n-\alpha+1}}\sum_{k=0}^{m-1} c_{\alpha,k}\left[ f^{(n-1)}( t_{m-k}) -f^{(n-1)}( t_{m-k-1} ) \right]+\mathcal{O}\left(
 dt^{n-\alpha+1} \right) \nonumber \\
 =&  \dfrac{dt^{n-\alpha-1}}{\gam{n-\alpha+1}}\left[f^{(n-1)}(t_m)-c_{\alpha,m-1}f^{(n-1)}(t_0)-  \sum_{k=1}^{m-1} \left( c_{\alpha,k-1}-c_{\alpha,k} \right) f^{(n-1)}( t_{m-k})\right]+\mathcal{O}\left(
 dt^{n-\alpha+1} \right).
\end{align}

It should be mentioned that the coefficients $c_{\alpha,k}$ of the previous expression are bounded and decreasing, which is exposed in the following proposition.

\begin{proposition}\label{prop:01}
The sequence $\set{c_{\alpha,k}}_{k=0}^\infty$, defined by \eqref{eq:S2-032}, is bounded and strictly decreasing for all $(n-1)<\alpha\leq n$.

\begin{proof}
To show that the sequence is bounded, we consider the following limit

\begin{eqnarray*}
\lim_{k\to \infty}\dfrac{(k+1)^{n-\alpha}}{k^{n-\alpha}}=\lim_{k\to \infty}\left(1+\dfrac{1}{k} \right)^{n-\alpha} \longrightarrow 1,
\end{eqnarray*}

as a consequence

\begin{eqnarray}
\lim_{k\to \infty}c_{\alpha,k}=\lim_{k\to \infty}\left[(k+1)^{n-\alpha}-k^{n-\alpha} \right]\longrightarrow 0.
\end{eqnarray}

On the other hand, to show that the sequence is strictly decreasing, we consider the following inequalities

\begin{eqnarray*}
k^{n-\alpha}<(k+1)^{n-\alpha}, & (k+1)^{n-\alpha}>0,
\end{eqnarray*}

then

\begin{eqnarray*}
\lim_{k\to \infty}\dfrac{(k+2)^{n-\alpha}+k^{n-\alpha}}{(k+1)^{n-\alpha}}<\lim_{k\to \infty} \dfrac{(k+2)^{n-\alpha}+(k+1)^{n-\alpha} }{(k+1)^{n-\alpha}} \longrightarrow 2,
\end{eqnarray*}

from the previous result

\begin{eqnarray*}
(k+2)^{n-\alpha}+k^{n-\alpha}<2(k+1)^{n-\alpha} & \Longrightarrow & (k+2)^{n-\alpha}-(k+1)^{n-\alpha}<(k+1)^{n-\alpha}-k^{n-\alpha},
\end{eqnarray*}

as a consequence

\begin{eqnarray}
\dfrac{c_{\alpha,k+1}}{c_{\alpha,k}}=\dfrac{(k+2)^{n-\alpha}-(k+1)^{n-\alpha}}{(k+1)^{n-\alpha}-k^{n-\alpha}}<1.
\end{eqnarray}

\end{proof}
\end{proposition}

Finally, from the equation \eqref{eq:S2-03} for the particular case $0<\alpha\leq 1$, we obtain the following expression

\begin{eqnarray}\label{eq:S2-04}
\ifr{C}{a}{D}{t}{\alpha}f(t)=\dfrac{dt^{-\alpha}}{\gam{2-\alpha}}\left[f(t_m)-c_{\alpha,m-1}f(t_0)-  \sum_{k=1}^{m-1} \left( c_{\alpha,k-1}-c_{\alpha,k} \right) f( t_{m-k})\right]+\mathcal{O}\left(
 dt^{2-\alpha} \right).
\end{eqnarray}

\section{Meshless Methods}

The meshless methods were created with the goal of eliminating some of the difficulties associated with constructing a mesh to generate a numerical approximation. In meshless methods, the approximation is built only from the nodes and this generates a computational time saving, since no time is wasted creating a mesh suitable for the problem we are trying to solve. One of the first meshless method was the Smoothed Particle Hydrodynamics Method \cite{lucy1977numerical,gingold1977smoothed}, designed to solve problems in astrophysics and, later, in fluid dynamics.

\subsection{Interpolation with Radial Basis Functions }

Let $\set{(x_j,u_j)}_{j=1}^{N_p}$ be a set of values, where $(x_j,u_j)\in \Omega \times \nset{R} \hspace{0.1cm} \forall j\geq 1$ with $\Omega \subset \nset{R}^d$. The interpolation problem in meshless methods is about finding a continuous function $\sigma: \Omega \subset \nset{R}^d \to \nset{R}$, such that

\begin{eqnarray}\label{eq:S2-01}
\sigma(x_j)=u_j , & \forall j\in \set{1,2,\cdots, N_p}.
\end{eqnarray}

In general, for the interpolation problem a function $\sigma$ is proposed as a linear combination using constants to be determined $\lambda_j\in \nset{R}$ and known base functions $B_j:\Omega \subset \nset{R}^d \to \nset{R}$, that is

\begin{eqnarray*}
\sigma(x)=\sum_{j=1}^{N_p}\lambda_j B_j(x), 
\end{eqnarray*}

then, from the interpolation condition \eqref{eq:S2-01}, the following matrix system is obtained

\begin{eqnarray}\label{eq:S2-02}
\begin{pmatrix}
B_1(x_1)& B_2(x_1)&\cdots &B_{N_p}(x_1)\\
B_1(x_2)&B_2(x_2)&\cdots & B_{N_p}(x_2)\\
\vdots & \vdots & \ddots & \vdots\\
B_1(x_{N_p})& B_2(x_{N_p})&\cdots & B_{N_p}(x_{N_p})
\end{pmatrix}
\begin{pmatrix}
\lambda_1\\ \lambda_2 \\ \vdots \\ \lambda_{N_p}
\end{pmatrix}=
\begin{pmatrix}
u_1\\
u_2\\
\vdots \\
u_{N_p}
\end{pmatrix},
\end{eqnarray}

which may be written in compact form as

\begin{eqnarray*}
G \Lambda =U,
\end{eqnarray*}

where $G_{jk}=B_k(x_j)$, $\Lambda_j=\lambda_j$ and $U_j=u_j$. It is said that the interpolation problem \eqref{eq:S2-02} is well posed, that is, the solution to the problem exists and is unique, if and only if the matrix $G$ is non-singular. 

The base functions $B_j$ are generally polynomial and trigonometric functions, which are computationally expensive to deal with larger-dimensional problems due to their dependence on geometric complexity. On the other hand, radial basis functions are constructed in terms of a distance, which makes them independent of the dimension of the problems, which gives them a clear advantage over other base functions. Before continuing it is necessary to have the following definition

\begin{definition}
Let $\Phi:\nset{R}^d \to \nset{R}$ be a function. Then, $\Phi$ is called radial, if there exists a function $\phi:\nset{R}_{\geq 0}\to \nset{R}$, such that

\begin{eqnarray*}
\Phi(x)=\phi(\norm{x}),
\end{eqnarray*}

where $\norm{ \ \cdot \ }:\nset{R}^d \to \nset{R}$ denotes any vector norm (generally the Euclidean norm).

\end{definition}

Let $\set{x_j}_{j=1}^{N_p}$ be a set of (random) nodes, then it is possible to construct a set of radial functions $\set{\Phi(x,x_j)}_{j=1}^{N_p}$, with

\begin{eqnarray*}
\Phi(x,x_j)=\phi\left(\norm{x-x_j}\right),
\end{eqnarray*}

therefore it is possible to generate a radial interpolant to implement the condition \eqref{eq:S2-01} as follows

\begin{eqnarray}\label{eq:S2-05}
\sigma(x)=\sum_{j=1}^{N_p}\lambda_j \Phi(x,x_j).
\end{eqnarray}

The methodology based on radial basis functions, proposed by Hardy \cite{hardy1971multiquadric}, arises from the need to apply multivariate interpolation in cartography problems using randomly dispersed nodes. Later, Kansa \cite{kansa1990multiquadrics1,kansa1990multiquadrics2} proposed to consider the analytical derivatives of radial basis functions to develop numerical schemes to solve partial differential equations.

\subsection{Solution of Differential Equations with Radial Basis Functions}

In this section we will give a brief introduction of how the radial basis functions methodology is used to solve a fractional partial differential equation, in the references \cite{gonzalez2016metodos,wendland,martinez2017applications1, martinez2019numerical,torreshern2019proposal}, it is possible to find more information and references to deepen the subject. Consider the following partial differential equation

\begin{eqnarray}\label{eq:S3-05}
\left\{
\begin{array}{cc}
\ifr{C}{0}{D}{t}{\alpha}u(x,t)-\mathcal{L}_{\beta,r}u(x,t)=u_I(x,t), & (x,t)\in \Omega \times D\\
u(x,t)=u_B(x,t), & (x,t)\in \partial \Omega \times D\\
u(x,t_0)=u_0(x), & x\in \Omega
\end{array}\right. ,
\end{eqnarray}

where the subscripts $I$ and $B$ refer to the interior and the border of the domain respectively. For the moment we focus on the fractional differential operator at interior of domain:

\begin{eqnarray*}
\ifr{C}{0}{D}{t}{\alpha}u(x,t)-\mathcal{L}_{\beta,r}u(x,t)=u_I(x,t),
\end{eqnarray*}

using the following notation

\begin{eqnarray*}
\left\{
\begin{array}{c}
\delta_{\alpha}:=\dfrac{\delta t^{-\alpha}}{\gam{2-\alpha}}\vspace{0.1cm} \\
u_I^m(x):=u_I(x,t_m) \vspace{0.1cm}\\
\mathcal{O}_{\alpha}(x):=\mathcal{O}\left(x, dt^{2-\alpha}\right)
\end{array}\right.,
\end{eqnarray*}

and considering \eqref{eq:S2-04}, we obtain that

\begin{eqnarray*}
\delta_{\alpha} \left[u^m(x)-c_{\alpha,m-1}u^0(x)-\sum_{k=1}^{m-1}\left( c_{\alpha,k-1}-c_{\alpha,k}  \right) u^{m-k}(x)   \right]-\mathcal{L}_{\beta,r}u^m(x)=u_I^m(x) +\mathcal{O}_{\alpha}^m(x),
\end{eqnarray*}

assuming $m\geq 1$, the previous expression we can write it as follows

\begin{eqnarray}\label{eq:S3-01}
\left( \delta_{\alpha}  -\mathcal{L}_{\beta,r} \right)  u^m(x)=u_I^m(x)+\delta_{\alpha}\left[ c_{\alpha,m-1}u^0(x)+\left( 1-\delta_{m-1,0} \right) \sum_{k=1}^{m-1}\left( c_{\alpha,k-1}-c_{\alpha,k}  \right) u^{m-k}(x)      \right]+\mathcal{O}_{\alpha}^m(x),
\end{eqnarray}

with $\delta_{m-1,0}$ the Kronecker delta and $u_I^m(x)=u_I(x,t_m)$. The superscript in $\mathcal{O}_{\alpha}^m$ is to indicate that it is the associated error of the approximation \eqref{eq:S2-04} to the time step $m$. As a consequence of the memory phenomenon of the fractional operator in time

\begin{eqnarray}\label{eq:S3-02}
\mathcal{O}_{\alpha}^m(x)=\mathcal{O}_{\alpha}^m\left(x,\mathcal{O}_{\alpha}^{m-1}(x),\mathcal{O}_{\alpha}^{m-2}(x),\cdots, \mathcal{O}_{\alpha}^1(x) \right),
\end{eqnarray}

so it is necessary to be careful with the value chosen for $m$, a very high value (that is, $0<dt\ll 1$) could lead to an error with an order of magnitude greater than expected. Once the equation \eqref{eq:S3-01} is obtained, it is necessary to define the conditions from which the values $u^m(x)$ are bounded, with which it is possible to determine its stability and convergence,  as shown in the references \cite{golbabai2019numerical,sun2006fully}. Before continuing,  we need to consider the following multi-index notation. Let $\nset{N}_0$ be the set $\nset{N}\cup\set{0}$, if $\gamma \in \nset{N}_0^d$ and $x\in  \nset{R}^d$, then

\begin{eqnarray*}
\left\{
\begin{array}{l}
 \abs{\gamma}:= \ds \sum_{k=1}^d \gamma_k\vspace{0.1cm}\\
\der{\partial}{x}{\gamma}:= \dfrac{\partial^{\abs{\gamma}}}{\partial x_1^{\gamma_1}\partial x_2^{\gamma_2}\cdots \partial x_d^{\gamma_d} }
\end{array}\right. ,
\end{eqnarray*}

considering  $\Omega \subset \nset{R}^d$ and using the previous notation, it is possible to define the following set of functions

\begin{eqnarray}
H^s(\Omega):=\set{  f(x)\in C^s(\Omega) \ : \  \der{\partial}{x}{\gamma}f(x) \in L^2(\Omega)  \hspace{0.1cm} \forall \abs{\gamma}\leq s },
\end{eqnarray}

it should be noted that in general, if $0<\beta \leq 1$, it is fulfills that

\begin{eqnarray}
\lim_{\beta \to 1} \mathcal{L}_{\beta,r}f(x) \longrightarrow \mathcal{L}_{1,r}f(x),
\end{eqnarray}

then if $f(x)\in H^2(\Omega)$, there exists $c>0$ such that

\begin{eqnarray}
\norm{\mathcal{L}_{\beta,r}f(x) }\leq c\norm{\mathcal{L}_{1,r}f(x) },
\end{eqnarray}

considering the above it is possible to prove the following proposition

\begin{proposition}
Let $\set{u^j(x)}_{j=1}^m$ be a sequence, defined by \eqref{eq:S3-01} on a domain $\Omega \subset \nset{R}^d$, with $u^j(x)\in H^2(\Omega) \ \forall j\geq 1$. Then for all $0<\alpha,\beta\leq 1$, it is fulfills that

\begin{eqnarray}\label{eq:S4-04}
\norm{\delta_{\alpha}   u^j(x)} \leq \dfrac{M}{c_{\alpha,j-1}}+ \norm{\delta_\alpha u^0(x)}, & j=1,2,\cdots, m,
\end{eqnarray}

where

\begin{eqnarray*}
M=\max_{1\leq k \leq m}\set{ \norm{u_I^k(x)} +\norm{\mathcal{L}_{\beta,r}u^k(x) }+\norm{\mathcal{O}_{\alpha}^k(x)}}.
\end{eqnarray*}

\begin{proof}

We proceed to prove \eqref{eq:S4-04} by induction:

\begin{itemize}

\item[i)]For the case $j=1$, from \eqref{eq:S3-01} we have that

\begin{eqnarray*}
\left( \delta_{\alpha}  -\mathcal{L}_{\beta,r} \right)  u^1(x)=u_I^1(x)+\delta_{\alpha}u^0(x) +\mathcal{O}_{\alpha}^1(x),
\end{eqnarray*}

then

\begin{eqnarray}\label{eq:S5-01}
\norm{\left( \delta_{\alpha}  -\mathcal{L}_{\beta,r} \right)  u^1(x)}\leq \norm{u_I^1(x)}+\norm{\delta_{\alpha}u^0(x)} +\norm{\mathcal{O}_{\alpha}^1(x)},
\end{eqnarray}

on the other hand, considering that $u^1(x)\in H^2(\Omega)$

\begin{eqnarray*}
\delta_{\alpha}u^1(x)=\left( \delta_{\alpha}  -\mathcal{L}_{\beta,r} \right)  u^1(x)+ \mathcal{L}_{\beta,r}u^1(x),
\end{eqnarray*}

then

\begin{eqnarray}\label{eq:S5-02}
\norm{\delta_{\alpha}u^1(x)}\leq \norm{\left( \delta_{\alpha}  -\mathcal{L}_{\beta,r} \right)  u^1(x)}+ \norm{\mathcal{L}_{\beta,r}u^1(x)},
\end{eqnarray}

as a consequence of \eqref{eq:S5-01} and \eqref{eq:S5-02}, we obtain that

\begin{eqnarray*}
\norm{\delta_{\alpha}u^1(x)}\leq \norm{u_I^1(x)}+ \norm{\mathcal{L}_{\beta,r}u^1(x)}+\norm{\mathcal{O}_{\alpha}^1(x)}+\norm{\delta_{\alpha}u^0(x)} ,
\end{eqnarray*}

therefore

\begin{eqnarray}
\norm{\delta_{\alpha}   u^1(x)} \leq \dfrac{M}{c_{\alpha,0}} +\norm{\delta_{\alpha}u^0(x)}.
\end{eqnarray}

\item[ii)] For the case $2\leq j \leq m-1$, we assume by induction hypothesis that it is fulfills that

\begin{eqnarray}\label{eq:S5-03}
\norm{\delta_{\alpha}   u^j(x)} \leq \dfrac{M}{c_{\alpha,j-1}} +\norm{\delta_{\alpha}u^0(x)}. 
\end{eqnarray}

\item[iii)] For the case $j=m$, from \eqref{eq:S3-01} we have that

\begin{eqnarray*}
\left( \delta_{\alpha}  -\mathcal{L}_{\beta,r} \right)  u^m(x)=u_I^m(x)+\delta_{\alpha}\left[ c_{\alpha,m-1}u^0(x)+\sum_{k=1}^{m-1}\left( c_{\alpha,k-1}-c_{\alpha,k}  \right) u^{m-k}(x)      \right]+\mathcal{O}_{\alpha}^m(x),
\end{eqnarray*}

in addition to the \textbf{Proposition \ref{prop:01}}, we have that $0<c_{k+1}<c_{k}$ if $ 0\leq k< \infty$, then

\begin{eqnarray}\label{eq:S5-04}
\norm{\left( \delta_{\alpha}  -\mathcal{L}_{\beta,r} \right)  u^m(x)}\leq \norm{u_I^m(x)}+c_{\alpha,m-1}\norm{\delta_{\alpha}u^0(x)}+\sum_{k=1}^{m-1}\left( c_{\alpha,k-1}-c_{\alpha,k}  \right) \norm{\delta_{\alpha}u^{m-k}(x)}    +\norm{\mathcal{O}_{\alpha}^m(x)},
\end{eqnarray}

on the other hand, considering that $u^m(x)\in H^2(\Omega)$

\begin{eqnarray*}
\delta_{\alpha}u^m(x)=\left( \delta_{\alpha}  -\mathcal{L}_{\beta,r} \right)  u^m(x)+ \mathcal{L}_{\beta,r}u^m(x),
\end{eqnarray*}

then

\begin{eqnarray}\label{eq:S5-05}
\norm{\delta_{\alpha}u^m(x)}\leq \norm{\left( \delta_{\alpha}  -\mathcal{L}_{\beta,r} \right)  u^m(x)}+ \norm{\mathcal{L}_{\beta,r}u^m(x)},
\end{eqnarray}

as a consequence of \eqref{eq:S5-04} and \eqref{eq:S5-05}, we obtain that

\begin{eqnarray*}
\norm{\delta_{\alpha}   u^m(x)} \leq\norm{ u_I^m(x)}+\norm{\mathcal{L}_{\beta,r}u^m(x)}+\norm{\mathcal{O}_{\alpha}^m(x)} +c_{\alpha,m-1}\norm{\delta_\alpha u^0(x)} +\sum_{k=1}^{m-1}\left(c_{\alpha,k-1}-c_{\alpha,k} \right)\norm{\delta_\alpha u^{m-k}(x)},
\end{eqnarray*}

then

\begin{eqnarray*}
\norm{\delta_{\alpha}   u^m(x)} \leq M +c_{\alpha,m-1}\norm{\delta_\alpha u^0(x)} +\sum_{k=1}^{m-1}\left(c_{\alpha,k-1}-c_{\alpha,k} \right)\norm{\delta_\alpha u^{m-k}(x)},
\end{eqnarray*}

as a consequence of the induction hypothesis \eqref{eq:S5-03}

\begin{eqnarray*}
\norm{\delta_{\alpha}   u^m(x)}  \leq M+c_{\alpha,m-1}\norm{\delta_\alpha u^0(x)}+\sum_{k=1}^{m-1}\left(c_{\alpha,k-1}-c_{\alpha,k} \right)\left( \dfrac{M}{c_{\alpha,m-k-1}} +\norm{\delta_{\alpha}u^0(x)}  \right),
\end{eqnarray*}

and from the \textbf{Proposition \ref{prop:01}}, we have that $0<c_{m-1}<c_{m-k-1}$ if $ 1 \leq k\leq m-1$, therefore

\begin{align}
\norm{\delta_{\alpha}   u^m(x)} \leq& M+c_{\alpha,m-1}\norm{\delta_\alpha u^0(x)}+\sum_{k=1}^{m-1}\left(c_{\alpha,k-1}-c_{\alpha,k} \right)\left( \dfrac{M}{c_{\alpha,m-1}} +\norm{\delta_{\alpha}u^0(x)}  \right) \nonumber\\
=& M+c_{\alpha,m-1}\norm{\delta_\alpha u^0(x)}+\left(c_{\alpha,0}-c_{\alpha,m-1} \right)\left( \dfrac{M}{c_{\alpha,m-1}} +\norm{\delta_{\alpha}u^0(x)}  \right)\nonumber\\
=& \dfrac{M}{c_{\alpha,m-1}} +\norm{\delta_{\alpha}u^0(x)} . 
\end{align}

\end{itemize}
\end{proof}

\end{proposition}

From the equation \eqref{eq:S3-01} and considering the boundary of the domain, we obtain the following system

\begin{eqnarray}\label{eq:S3-03}
\widetilde{\mathcal{L}}_{\alpha,\beta,r}u^m(x)=\widetilde{u}_{\alpha,IB}^m(x)+ \mathcal{O}_{\alpha,\Omega}^m(x),
\end{eqnarray}

where

\begin{eqnarray*}
\widetilde{\mathcal{L}}_{\alpha,\beta,r}u^m(x):=\left\{
\begin{array}{cc}
\left( \delta_{\alpha}  -\mathcal{L}_{\beta,r} \right)  u^m(x), & \mbox{ if }x\in \Omega\\
u^m(x), & \mbox{ if }x\in \partial \Omega
\end{array}
\right.,
\end{eqnarray*}

\begin{eqnarray*}
\widetilde{u}_{\alpha,IB}^m(x):=\left\{
\begin{array}{cc}
\ds u_I^m(x)+\delta_{\alpha}\left[ c_{\alpha,m-1}u^0(x)+\left( 1-\delta_{m-1,0} \right) \sum_{k=1}^{m-1}\left( c_{\alpha,k-1}-c_{\alpha,k}  \right) u^{m-k}(x)      \right], & \mbox{ if }x\in \Omega\\
u_B^m(x), & \mbox{ if }x\in \partial \Omega
\end{array}
\right.,
\end{eqnarray*}

\begin{eqnarray*}
\mathcal{O}_{\alpha ,\Omega}^m(x):=\left\{
\begin{array}{cc}
\mathcal{O}_{\alpha}^m(x), & \mbox{ if }x\in \Omega\\
0, & \mbox{ if }x\in \partial \Omega
\end{array}
\right..
\end{eqnarray*}

Now considering a radial interpolant

\begin{eqnarray*}
\sigma^m(x)=\sum_{j=1}^{N_p}\lambda_j^m \Phi(x,x_j),
\end{eqnarray*}

and a set of (random) nodes $\set{x_j}_{j=1}^{N_p}\subset \overline{\Omega} $. Then,  substituting the interpolant $\sigma^m_I$ in the equation \eqref{eq:S3-03}, for each value of $x_j$, an interpolation condition analogous to \eqref{eq:S2-01} is obtained. Therefore we obtain the following matrix system

\begin{eqnarray}\label{eq:S3-04}
\begin{pmatrix}
\widetilde{\mathcal{L}}_{\alpha,\beta,r}\Phi_{11}&\widetilde{\mathcal{L}}_{\alpha,\beta,r} \Phi_{12}&\cdots &\widetilde{\mathcal{L}}_{\alpha,\beta,r}\Phi_{1N_p}\\
\widetilde{\mathcal{L}}_{\alpha,\beta,r}\Phi_{21}& \widetilde{\mathcal{L}}_{\alpha,\beta,r}\Phi_{22}&\cdots &\widetilde{\mathcal{L}}_{\alpha,\beta,r}\Phi_{2N_p}\\
\vdots & \vdots & \ddots & \vdots\\
\widetilde{\mathcal{L}}_{\alpha,\beta,r}\Phi_{N_p1}&\widetilde{\mathcal{L}}_{\alpha,\beta,r} \Phi_{N_p2}&\cdots &\widetilde{\mathcal{L}}_{\alpha,\beta,r}\Phi_{N_pN_p}\\
\end{pmatrix}
\begin{pmatrix}
\lambda_1^m\\ \lambda_2^m \\ \vdots \\ \lambda_{N_p}^m
\end{pmatrix}=
\begin{pmatrix}
\widetilde{u}_{\alpha,IB,1}^m+\mathcal{O}^m_{\alpha,\Omega,1}\\
\widetilde{u}_{\alpha,IB,2}^m+\mathcal{O}^m_{\alpha,\Omega,2}\\
\vdots \\
\widetilde{u}_{\alpha,IB,N_p}^m+\mathcal{O}^m_{\alpha,\Omega,N_p}
\end{pmatrix},
\end{eqnarray}

where

\begin{eqnarray*}
\left\{
\begin{array}{c}
\widetilde{\mathcal{L}}_{\alpha,\beta,r}\Phi_{ij}=\widetilde{\mathcal{L}}_{\alpha,\beta,r}\Phi(x_i,x_j) \vspace{0.1cm} \\
\widetilde{u}_{\alpha,IB,j}^m=\widetilde{u}_{\alpha,IB}^m(x_j) \vspace{0.1cm}\\
\mathcal{O}^m_{\alpha,\Omega,j}=\mathcal{O}^m_{\alpha,\Omega}(x_j)
\end{array}\right..
\end{eqnarray*}

Under the assumption that the above matrix is invertible, the interpolant may be written as

\begin{eqnarray*}
\sigma^m(x)=\sum_{j=1}^{N_p}\lambda_j^m \Phi(x,x_j)=\sum_{j=1}^{N_p}\left(\widetilde{ \lambda}_j^m +\widetilde{\mathcal{O}}_{\alpha,\Omega,j}^m \right) \Phi(x,x_j).
\end{eqnarray*}

From the previous expression, it becomes clear that the number of nodes chosen to find the solution is also a factor in which care must be taken when considering the errors of the solution. Assuming that the system \eqref{eq:S3-05} has an analytical solution  $u_s(x,t)$, we have that

\begin{eqnarray*}
\norm{\sigma^m(x)-u_s(x,t_m)}\leq N_p\max_{1\leq j\leq N_p}\set{\abs{ \widetilde{\mathcal{O}}_{ \alpha,\Omega,j}^m }  \norm{ \Phi(x,x_j)} },
\end{eqnarray*}

where in general

\begin{eqnarray*}
\lim_{dt\to 0} \abs{ \widetilde{\mathcal{O}}_{ \alpha,\Omega,j}^m } \to 0.
\end{eqnarray*}

Considering that the system \eqref{eq:S3-05} for $0<\alpha,\beta< 1$, in general has no analytical solution, we will use the root mean squared error of the operator $\widetilde{\mathcal{L}}_{\alpha,\beta,r}$ applied to the interpolant $\sigma^m(x)$  with the interpolation condition $\widetilde{u}_{\alpha,IB}(x_j)$ to estimate the error of the solution, that is,

\begin{eqnarray}
RMSE_m=\sqrt{\dfrac{1}{N_P}\sum_{j=1}^{N_p}\left(\widetilde{\mathcal{L}}_{\alpha,\beta,r}\sigma_j^m-\widetilde{u}_{\alpha,IB,j}^m \right)^2  }.
\end{eqnarray}

The system \eqref{eq:S3-04} may be written compactly as follows

\begin{eqnarray*}
G_{\alpha,\beta}\Lambda^m=U_{\alpha}^m,
\end{eqnarray*}

it is necessary to mention that in general, the matrix $G_{\alpha,\beta}$ fulfills the following condition

\begin{eqnarray*}
\lim_{N_p\to \infty}\det(G_{\alpha,\beta})\to 0,
\end{eqnarray*}

as a consequence, although $\det(G_{\alpha,\beta})\neq 0$, there is a risk that the matrix  is analytically invertible but numerically singular. To solve this problem, a  preconditioning matrix  $P$ is generated through the factorization $QR$ of the matrix $G_{\alpha,\beta}$ \cite{stoer2013}, that is,

\begin{eqnarray*}
G_{\alpha,\beta}=QR,
\end{eqnarray*}

then the following matrix is defined

\begin{eqnarray*}
\widetilde{Q}:=(\widetilde{Q}_{ij})=\left( \log\left( \exp(Q_{ij})+\dfrac{1}{\cond(G_{\alpha,\beta})} \right) \right),
\end{eqnarray*}

and the system \eqref{eq:S3-04} is replaced by the following system

\begin{eqnarray}\label{eq:S3-06}
\widetilde{ G}_{\alpha,\beta}\Lambda^m=\widetilde{U}_{\alpha}^m,
\end{eqnarray}

where

\begin{align*}
\widetilde{ G}_{\alpha,\beta}:=&PG_{\alpha,\beta}=\left(\widetilde{Q}R \right)^{-1}G_{\alpha,\beta},\\
\widetilde{U}_{\alpha}^m:=&PU_{\alpha}^m=\left(\widetilde{Q}R \right)^{-1}U_{\alpha}^m,
\end{align*}

with which the following relationship between the matrices $G_{\alpha,\beta}$ and $\widetilde{ G}_{\alpha,\beta}$  is guaranteed

\begin{eqnarray*}
\dfrac{1}{\cond(G_{\alpha,\beta})}<\dfrac{\cond(\widetilde{ G}_{\alpha,\beta})}{\cond(G_{\alpha,\beta})} \ll 1 .
\end{eqnarray*}

\subsubsection{Examples}

For a set  of chosen (random) nodes $\set{x_j}_{j=1}^{N_p}$, a set of radial functions $\set{\Phi(x,x_j)}_{j=1}^{N_p}$ is generated, where

\begin{eqnarray}\label{eq:S4-03}
\Phi(x,x_j)=\norm{x-x_j}_2^3.
\end{eqnarray}

The following examples are solved using the set of radial functions above and the system \eqref{eq:S3-06}, with the following particular values

\begin{eqnarray*}
\begin{array}{ccc}
\widetilde{\sigma}=0.25, & \widetilde{r}=0.05, & dt=\dfrac{1}{25}.
\end{array}
\end{eqnarray*}

\begin{example}
\begin{eqnarray}\label{eq:S4-01}
\left\{
\begin{array}{cc}
\ifr{C}{0}{D}{t}{\alpha}u(x,t)-\mathcal{L}_{\beta,r}u(x,t)=u_I(x,t), & (x,t)\in [0,1] \times [0,1]\\
u(x,t)=0, & (x,t)\in \partial ( [0,1] ) \times [0,1]\\
u(x,0)=(1-x)\sin^2(x), & x\in [0,1]
\end{array}\right. ,
\end{eqnarray}

where

\begin{eqnarray*}
\mathcal{L}_{\beta,r}:= \dfrac{1}{2}\widetilde{\sigma}^2  \ifr{}{0}{D}{r}{\beta+1}+\left(\widetilde{r}-\dfrac{1}{2}\widetilde{\sigma}^2  \right) \ifr{}{0}{D}{r}{\beta}-\widetilde{r},
\end{eqnarray*}

with

\begin{align*}
u_I(x,t)=&\widetilde{\sigma}^2 \left[\sin(2 x)-(1 - x) \cos(2 x)    \right](t+1)^2 \\
&+\left[ 2(t+1)+\widetilde{r}(t+1)^2 \right](1-x)\sin^2(x)\\
&+\left(\widetilde{r}-\dfrac{1}{2}\widetilde{\sigma}^2  \right)  \left[\sin^2(x)- (1 - x) \sin(2x)\right](t+1)^2,
\end{align*}

and whose analytical solution for the particular case $\alpha=\beta=1$ is the following

\begin{eqnarray*}
u(x,t)=(t+1)^2(1-x)\sin^2(x).
\end{eqnarray*}

Different numbers of Chebyshev nodes are used to solve the system of equations \eqref{eq:S4-01} (see Figure \ref{fig:01}). The numerical solutions for different values of $\alpha$ and $\beta$ for $100$ Chevichev nodes are presented in Figure \ref{fig:02}, and some results are shown in Table \ref{tab:01}.

\begin{figure}[!ht]
\centering
\begin{subfigure}[c]{0.33\linewidth}
\includegraphics[width=\linewidth,height=0.46\linewidth]{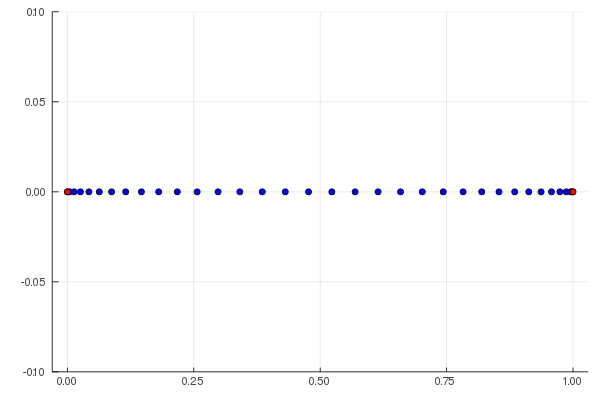}
\caption{$N_p=36$.}
\end{subfigure}
\begin{subfigure}[c]{0.33\linewidth}
\includegraphics[width=\linewidth,height=0.46\linewidth]{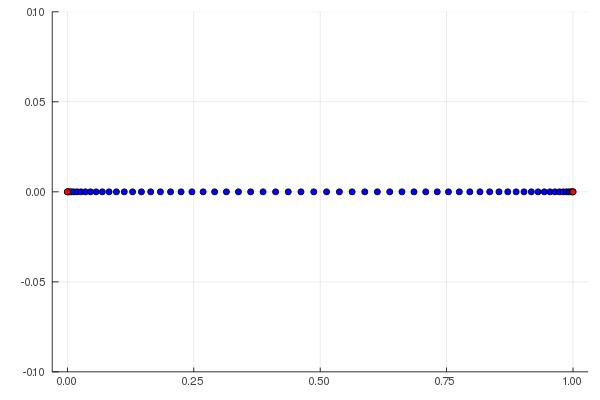}
\caption{$N_p=64$.}
\end{subfigure}
\begin{subfigure}[c]{0.33\linewidth}
\includegraphics[width=\linewidth,height=0.46\linewidth]{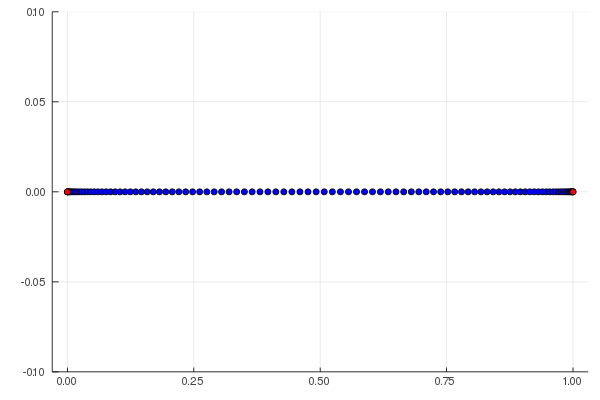}
\caption{$N_p=100$.}
\end{subfigure}
\caption{Different numbers of Chebyshev nodes used.}\label{fig:01}
\end{figure}

\begin{figure}[!ht]
\centering
\begin{subfigure}[c]{0.33\linewidth}
\includegraphics[width=\linewidth,height=0.65\linewidth]{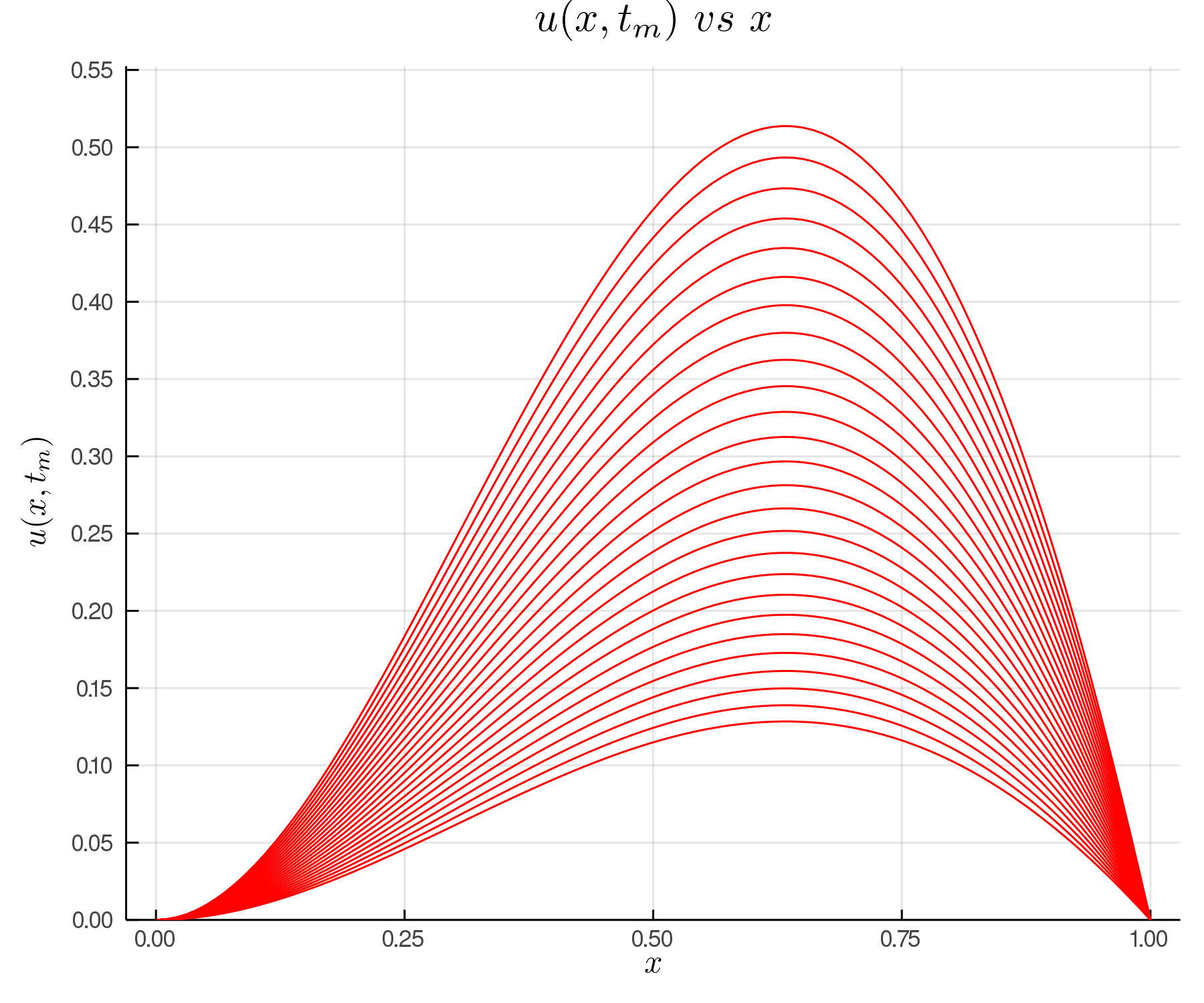}
\caption{\footnotesize Analytical solution.}
\end{subfigure}
\begin{subfigure}[c]{0.33\linewidth}
\includegraphics[width=\linewidth,height=0.65\linewidth]{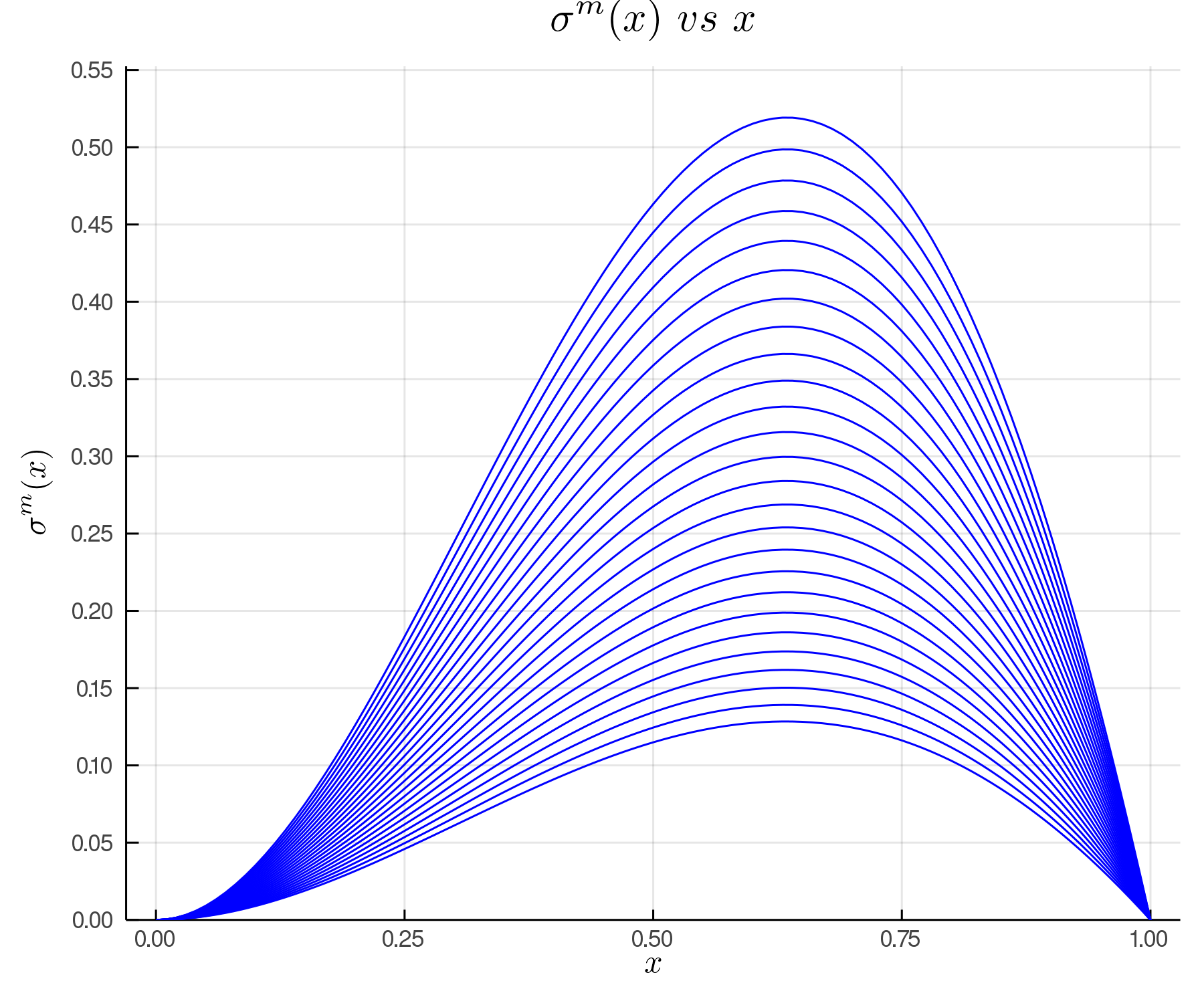}
\caption{\footnotesize Numerical solution for $(\alpha,\beta)=(1,1)$.}
\end{subfigure}
\begin{subfigure}[c]{0.33\linewidth}
\includegraphics[width=\linewidth,height=0.65\linewidth]{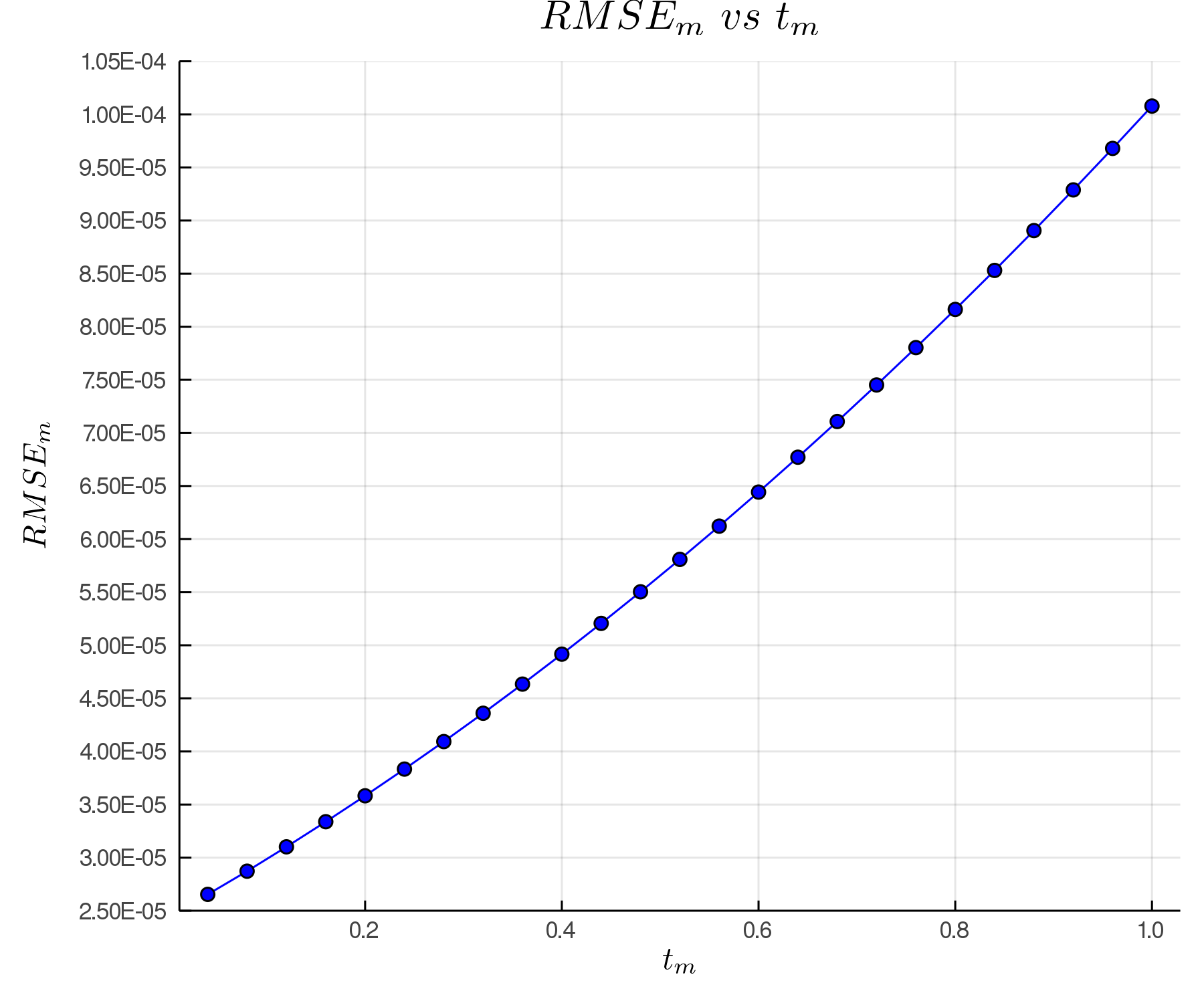}
\caption{\footnotesize $RMSE_m$ for $(\alpha,\beta)=(1,1)$.}
\end{subfigure}
\\
\begin{subfigure}[c]{0.33\linewidth}
\includegraphics[width=\linewidth,height=0.65\linewidth]{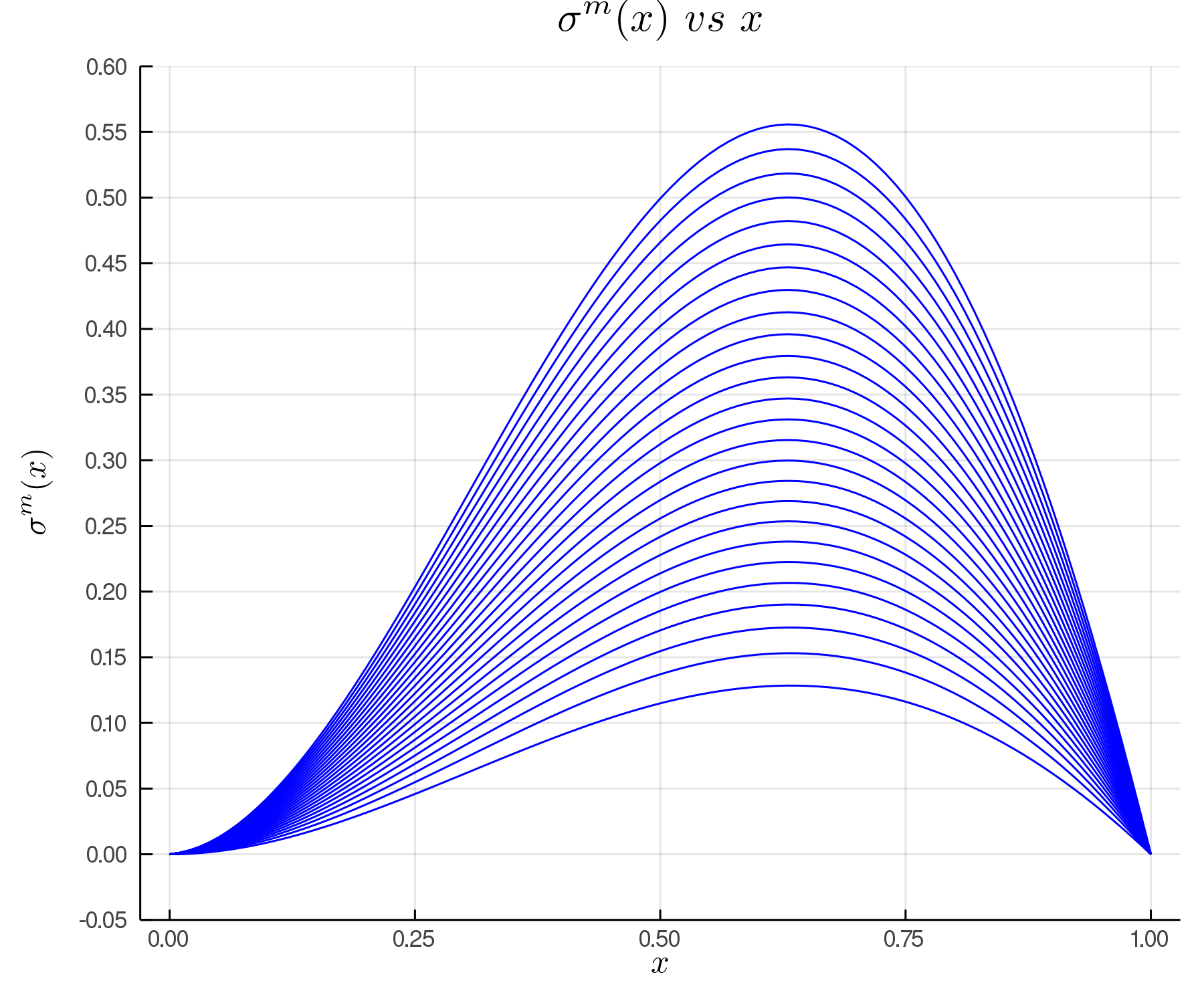}
\caption{\footnotesize Numerical solution for $(\alpha,\beta)=(0.7,1)$.}
\end{subfigure}
\begin{subfigure}[c]{0.33\linewidth}
\includegraphics[width=\linewidth,height=0.65\linewidth]{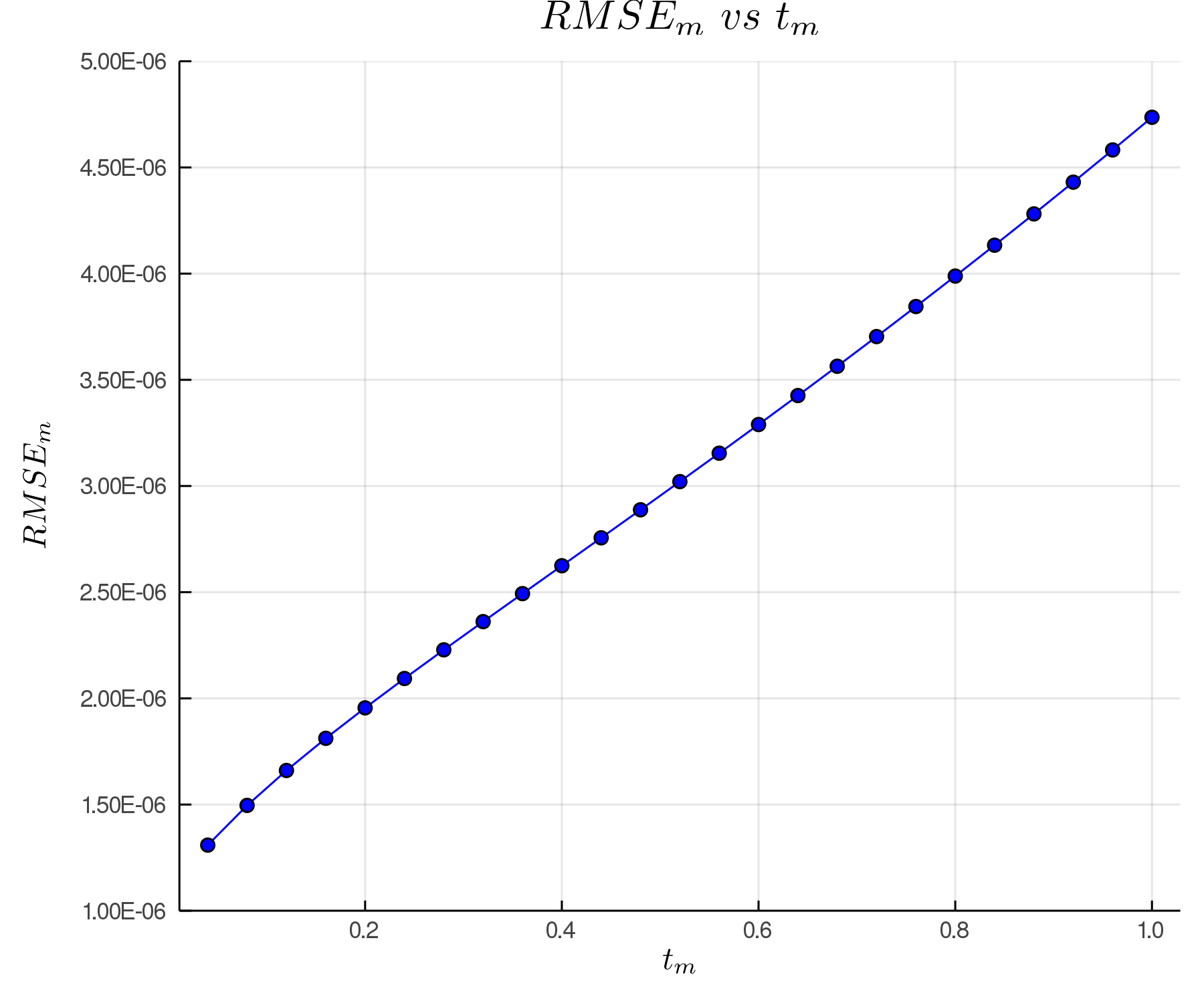}
\caption{\footnotesize $RMSE_m$ for $(\alpha,\beta)=(0.7,1)$.}
\end{subfigure}
\\
\begin{subfigure}[c]{0.33\linewidth}
\includegraphics[width=\linewidth,height=0.65\linewidth]{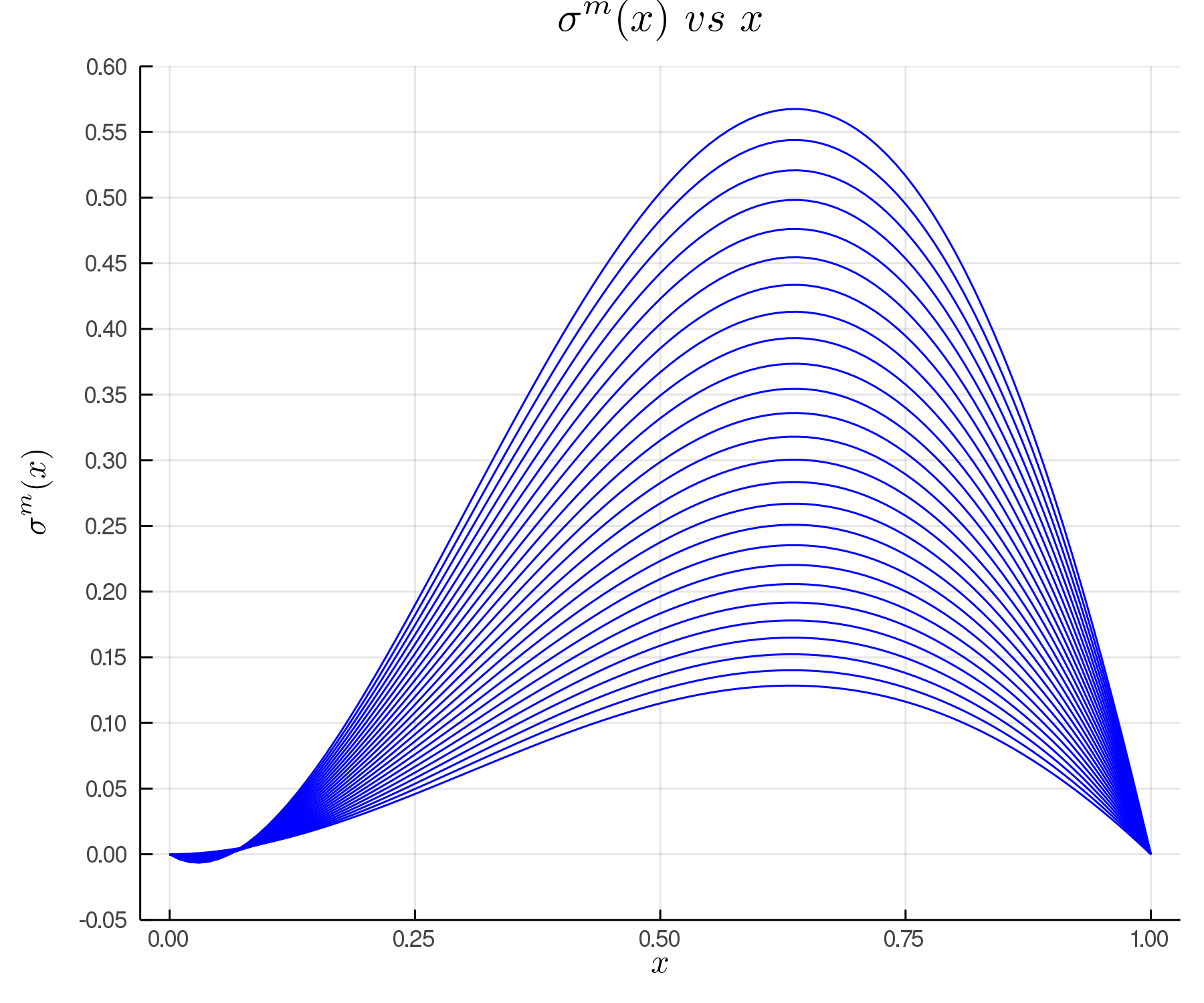}
\caption{\footnotesize Numerical solution for $(\alpha,\beta)=(1,0.75)$.}
\end{subfigure}
\begin{subfigure}[c]{0.33\linewidth}
\includegraphics[width=\linewidth,height=0.65\linewidth]{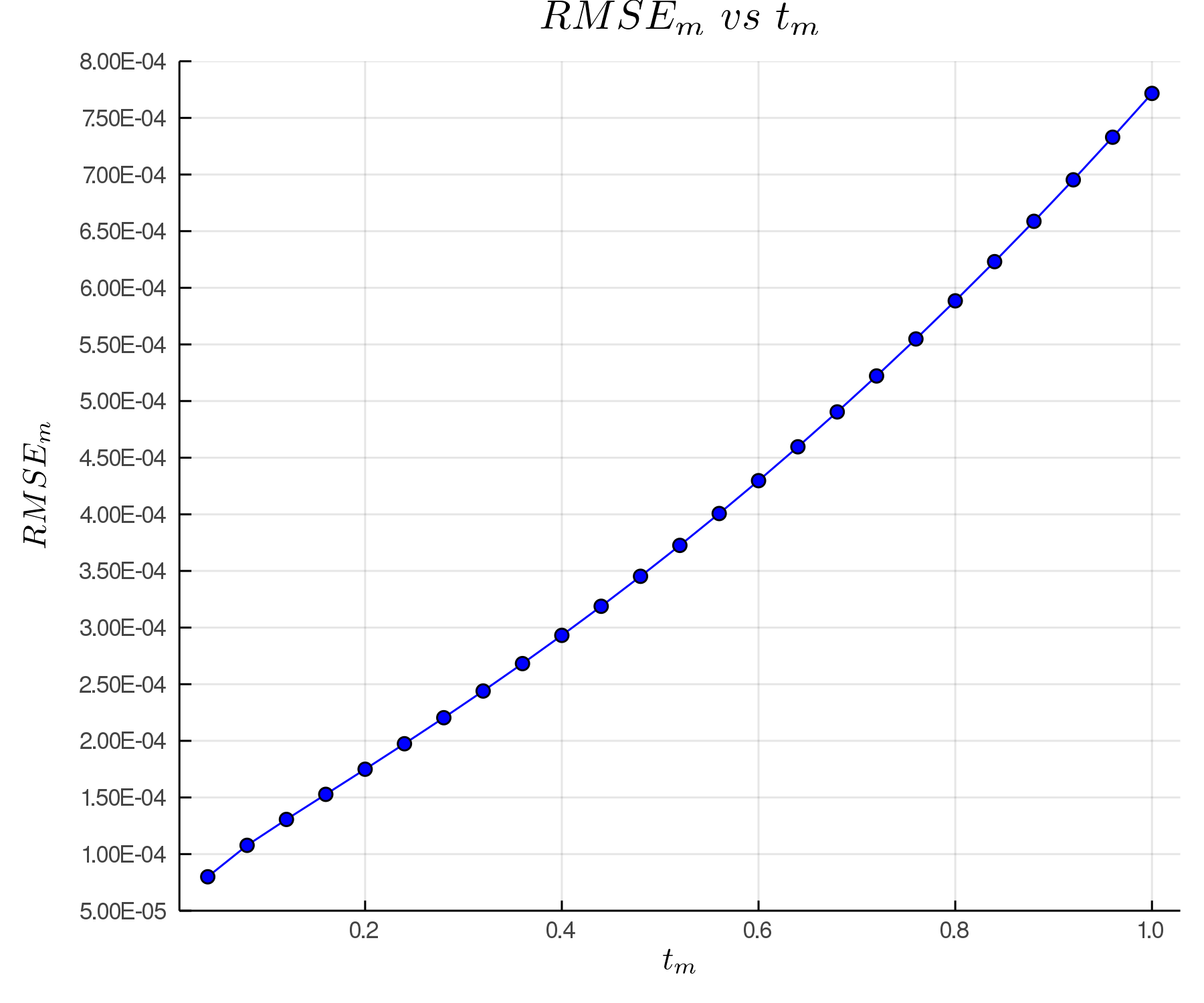}
\caption{\footnotesize $RMSE_m$ for $(\alpha,\beta)=(1,0.75)$.}
\end{subfigure}
\\
\begin{subfigure}[c]{0.33\linewidth}
\includegraphics[width=\linewidth,height=0.65\linewidth]{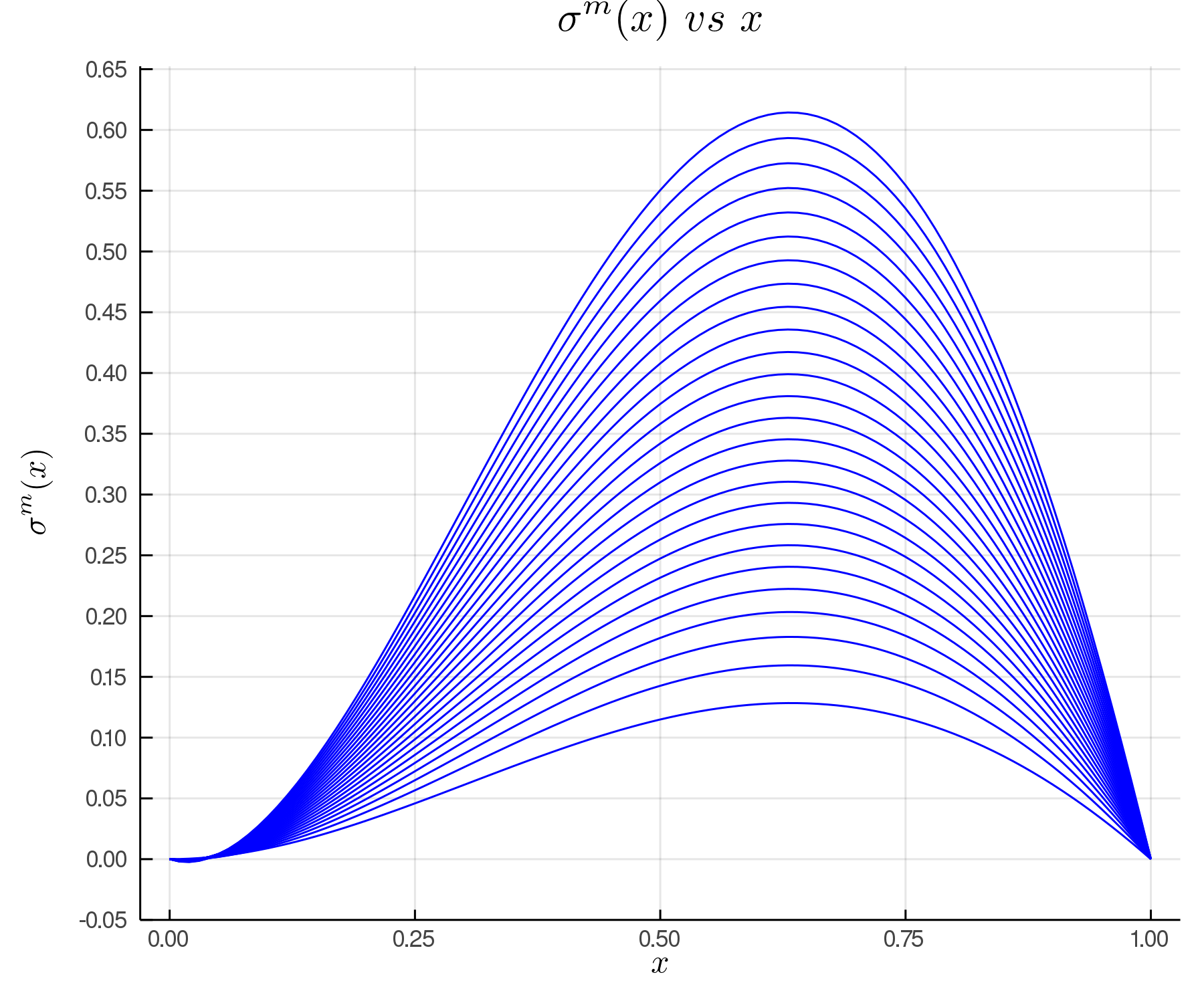}
\caption{ \footnotesize Numerical solution for $(\alpha,\beta)=(0.65,0.8)$.}
\end{subfigure}
\begin{subfigure}[c]{0.33\linewidth}
\includegraphics[width=\linewidth,height=0.65\linewidth]{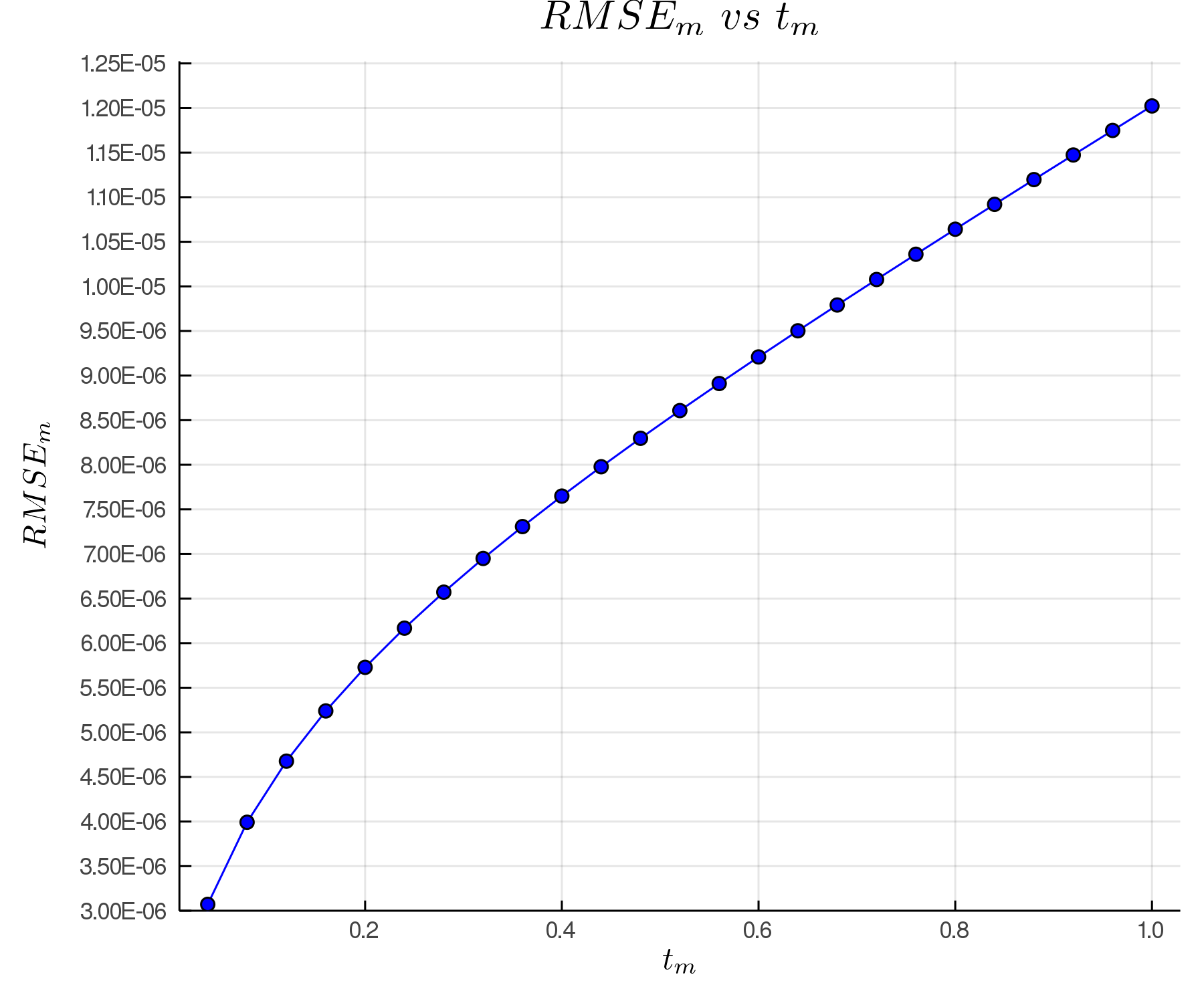}
\caption{\footnotesize $RMSE_m$ for $(\alpha,\beta)=(0.65,0.8)$.}
\end{subfigure}
\caption{The analytical solution and the numerical solutions with respect to space for different moments in time are presented. The $RMSE$ is presented with respect to time for the different numerical solutions.}\label{fig:02}
\end{figure}

\begin{table}[!ht]
\centering
$
\small
\begin{array}{cccccc}
\toprule
\alpha&\beta&N_p& \cond(G_{\alpha,\beta})&\cond(\widetilde{ G}_{\alpha,\beta})&RMSE\\
\midrule
          &       & 36    & 3.19985E+07 & 9.03438E+00 & 6.07019E-07 \\
    1     & 1     & 64    & 1.89809E+08 & 2.03053E+01 & 4.26112E-06 \\
          &       & 100   & 7.41018E+08 & 3.52912E+01 & 1.00781E-04 \\ \midrule
          &       & 36    & 5.42967E+06 & 5.34347E+00 & 9.65987E-08 \\
    0.7   & 1     & 64    & 3.21286E+07 & 1.37232E+01 & 6.32924E-07 \\
          &       & 100   & 1.25345E+08 & 2.48729E+01 & 4.73620E-06 \\ \midrule
          &       & 36    & 1.42516E+08 & 1.04670E+01 & 6.59344E-06 \\
    1     & 0.75  & 64    & 9.92520E+08 & 2.46329E+01 & 1.04166E-04 \\
          &       & 100   & 4.35794E+09 & 4.49581E+01 & 7.71606E-04 \\ \midrule
          &       & 36    & 1.24173E+07 & 5.31271E+00 & 5.86635E-07 \\
    0.65  & 0.8   & 64    & 8.34941E+07 & 1.17088E+01 & 6.29195E-06 \\
          &       & 100   & 3.58730E+08 & 2.53325E+01 & 1.20229E-05 \\ \bottomrule
\end{array}
$
\caption{Values obtained for the different numerical solutions, the value of $RMSE$ is presented for the final time step.}\label{tab:01}
\end{table}

\end{example}

\newpage

\begin{example}
\begin{eqnarray}\label{eq:S4-02}
\left\{
\begin{array}{cc}
\ifr{C}{0}{D}{t}{\alpha}u(x,y,t)-\mathcal{L}_{\beta,r}u(x,y,t)=u_I(x,y,t), & (x,y,t)\in [0,1]\times [0,1] \times [0,1]\\
u(x,y,t)=u_B(x,y,t), & (x,y,t)\in \partial( [0,1]\times [0,1]) \times [0,1]\\
u(x,y,0)=\dfrac{1}{4}(1-x^2-y^2)(2-x^2-y^2)\sin^2\left(2(x^2+y^2)\right), & (x,y)\in [0,1]\times [0,1]
\end{array}\right. ,
\end{eqnarray}

where

\begin{eqnarray*}
\mathcal{L}_{\beta,r}:= \dfrac{1}{2}\widetilde{\sigma}^2  \ifr{}{0}{D}{r}{\beta+1}+\left(\widetilde{r}-\dfrac{1}{2}\widetilde{\sigma}^2  \right) \ifr{}{0}{D}{r}{\beta}-\widetilde{r},
\end{eqnarray*}

with

\begin{eqnarray*}
\footnotesize
\begin{array}{c}
u_B(x,y,t)=\left\{
\begin{array}{cc}
\dfrac{1}{4}(t+1)^2(1-y^2)(2-y^2)\sin^2\left(2y^2\right), &\mbox{ if } (x,y,t)\in \set{0}\times [0,1]\times [0,1] \vspace{0.1cm}\\
\dfrac{1}{4}(t+1)^2y^2(y^2-1)\sin^2\left(2(1+y^2)\right), &\mbox{ if } (x,y,t)\in \set{1}\times [0,1]\times [0,1]\vspace{0.1cm}\\
\dfrac{1}{4}(t+1)^2(1-x^2)(2-x^2)\sin^2\left(2x^2\right), &\mbox{ if } (x,y,t)\in  [0,1]\times \set{0}\times [0,1]\vspace{0.1cm}\\
\dfrac{1}{4}(t+1)^2x^2(x^2-1)\sin^2\left(2(1+x^2)\right), &\mbox{ if } (x,y,t)\in [0,1]\times \set{1}\times [0,1]
\end{array}\right.,
\end{array}
\end{eqnarray*}

and

\begin{align*}
\footnotesize
\begin{array}{rl}
u_I(x,y,t)=&\dfrac{\widetilde{\sigma}^2}{8} \left\{ 3-\left[3 + 58 (x^2+y^2) - 96(x^2+y^2)^2 + 32 (x^2+y^2)^3\right] \cos\left(4 (x^2+y^2)\right)\right\}(t+1)^2 \vspace{0.1cm} \\
&-\dfrac{\widetilde{\sigma}^2 }{4}\left\{3 (x^2+y^2)+ \left[4 - 30 (x^2+y^2) + 18 (x^2+y^2)^2\right] \sin\left(4 (x^2+y^2)\right)     \right\}(t+1)^2 \vspace{0.1cm} \\
&+ \dfrac{1}{4}\left[ 2(t+1)+\widetilde{r}(t+1)^2 \right](1-x^2-y^2)(2-x^2-y^2)\sin^2(2(x^2+y^2)) \vspace{0.1cm} \\
&+\dfrac{1}{2}\left(\widetilde{r}-\dfrac{1}{2}\widetilde{\sigma}^2  \right) \sqrt{x^2+y^2} \left[3 - 2(x^2+y^2)\right]\sin^2\left(2 (x^2+y^2)\right)(t+1)^2 \vspace{0.1cm}\\
&-\left(\widetilde{r}-\dfrac{1}{2}\widetilde{\sigma}^2  \right)  \sqrt{x^2+y^2}  \left[2 - 3 (x^2+y^2) + (x^2+y^2)^2\right] \sin\left(4 (x^2+y^2)\right) (t+1)^2,
\end{array}
\end{align*}

whose analytical solution for the particular case $\alpha=\beta=1$ is the following

\begin{eqnarray*}
u(x,y,t)=\frac{1}{4}(t+1)^2(1-x^2-y^2)(2-x^2-y^2)\sin^2\left(2(x^2+y^2)\right).
\end{eqnarray*}

For this example, we use a combination of Halton type nodes within the domain and Cartesian nodes at the boundary. Different numbers of nodes are used to solve the system of equations \eqref{eq:S4-02} (see Figure \ref{fig:03}). The numerical solutions for different values of $\alpha$ and $\beta$ for $400$ nodes are presented in Figure \ref{fig:04}, and some results are shown in Table \ref{tab:02}.

\begin{figure}[!ht]
\centering
\begin{subfigure}[c]{0.33\linewidth}
\includegraphics[width=\linewidth,height=0.84\linewidth]{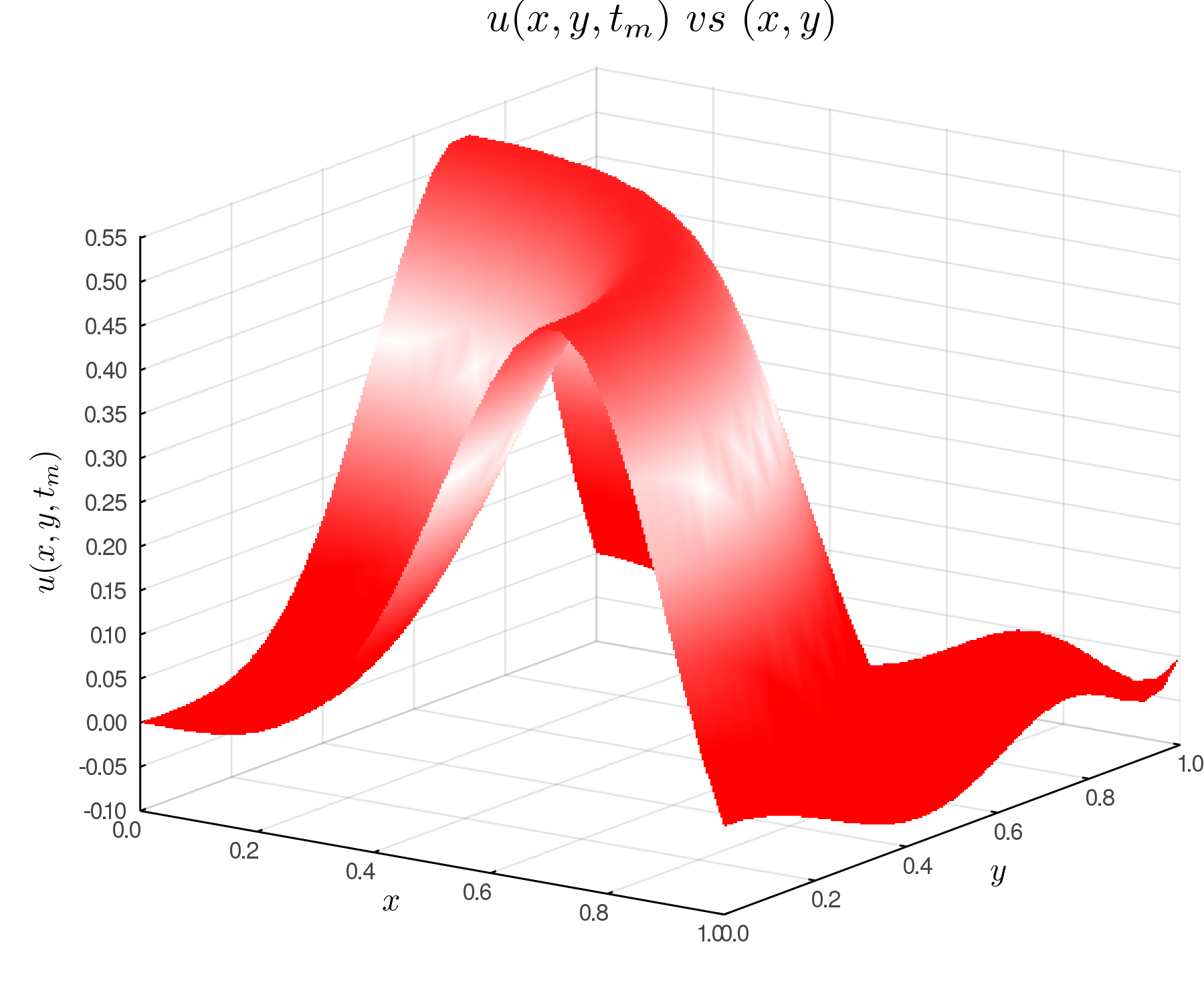}
\caption{\footnotesize Analytical solution.}
\end{subfigure}
\begin{subfigure}[c]{0.33\linewidth}
\includegraphics[width=\linewidth,height=0.84\linewidth]{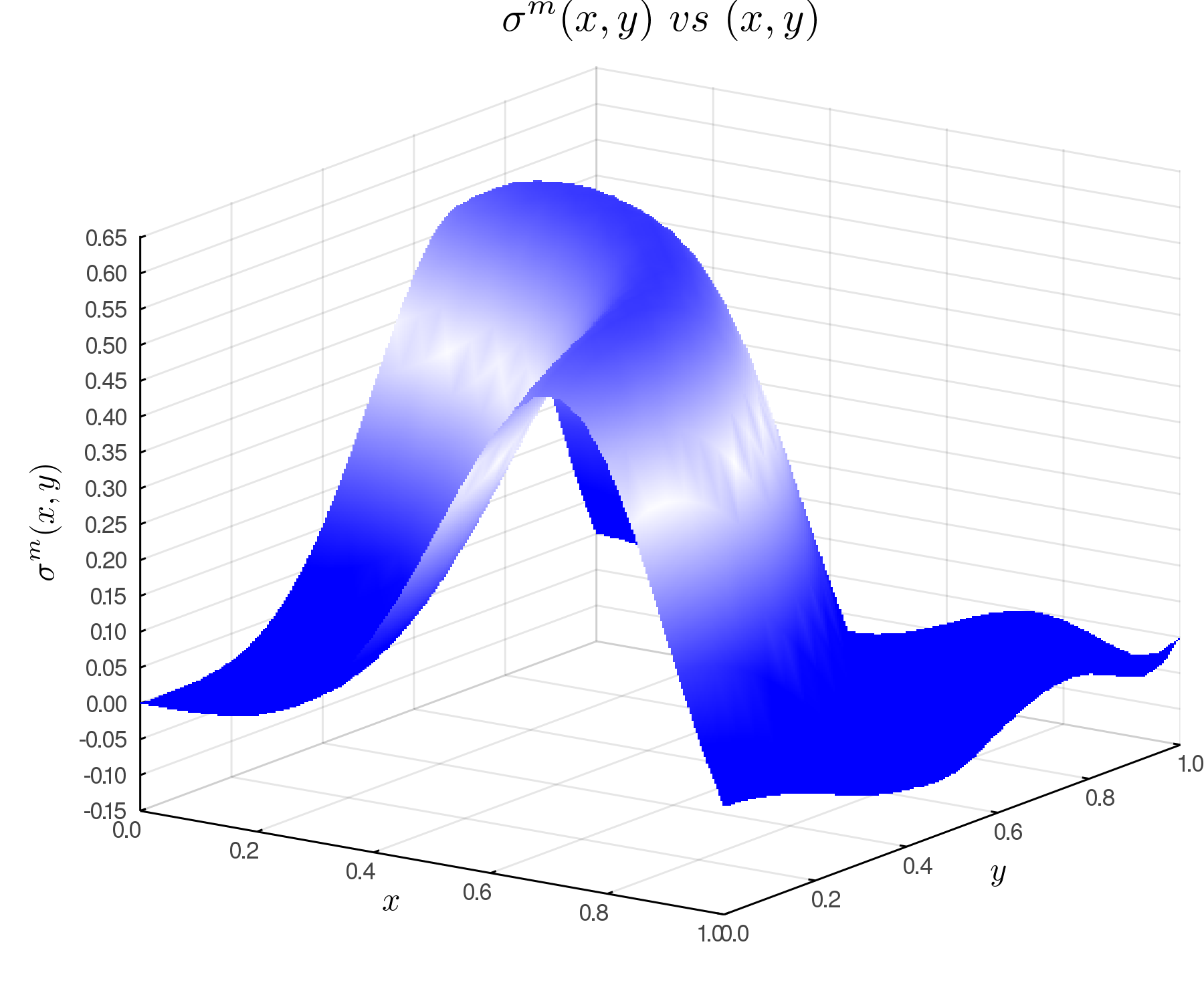}
\caption{\footnotesize Numerical solution for $(\alpha,\beta)=(1,1)$.}
\end{subfigure}
\begin{subfigure}[c]{0.33\linewidth}
\includegraphics[width=\linewidth,height=0.84\linewidth]{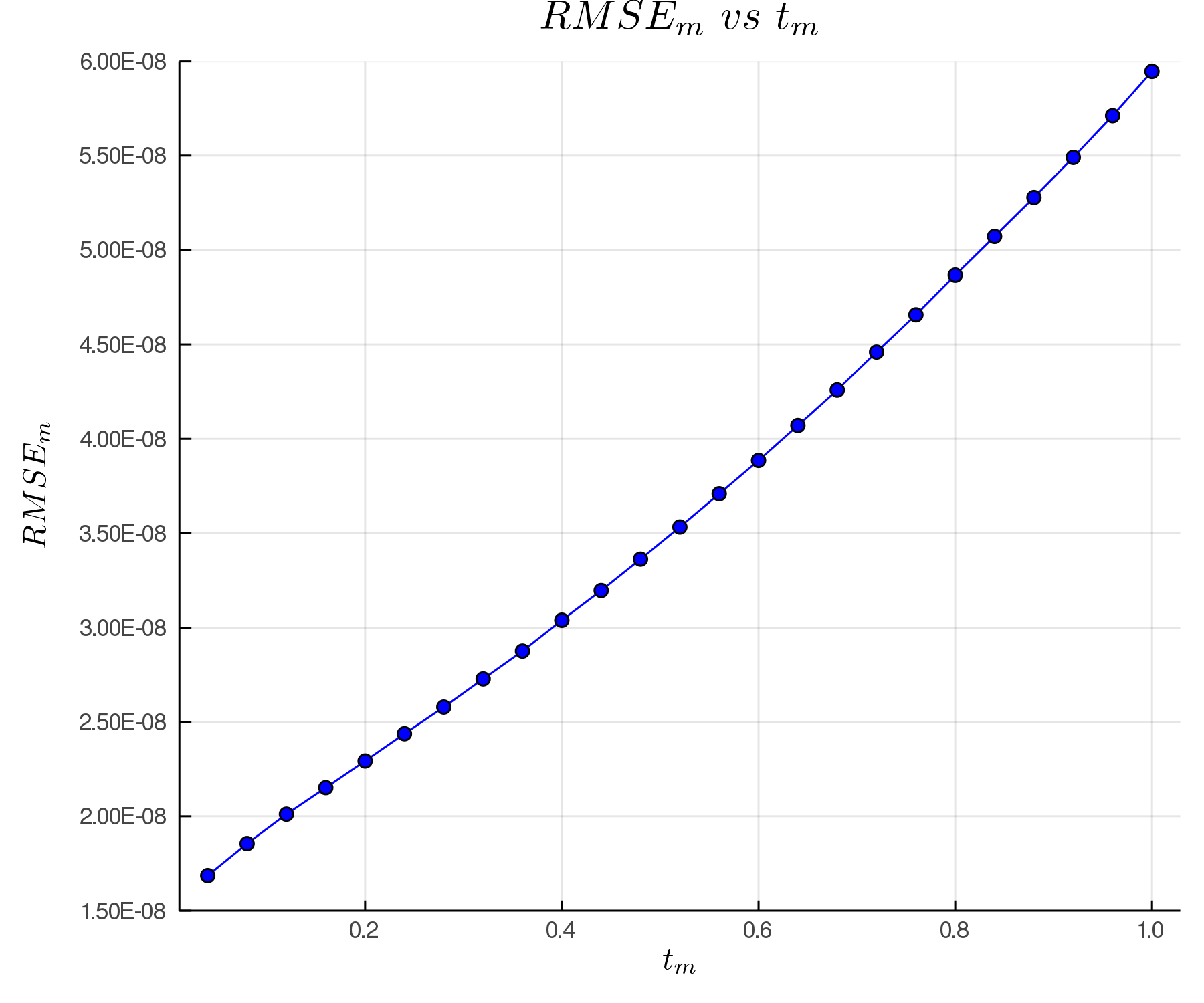}
\caption{\footnotesize $RMSE_m$ for $(\alpha,\beta)=(1,1)$.}
\end{subfigure}
\\
\begin{subfigure}[c]{0.33\linewidth}
\includegraphics[width=\linewidth,height=0.84\linewidth]{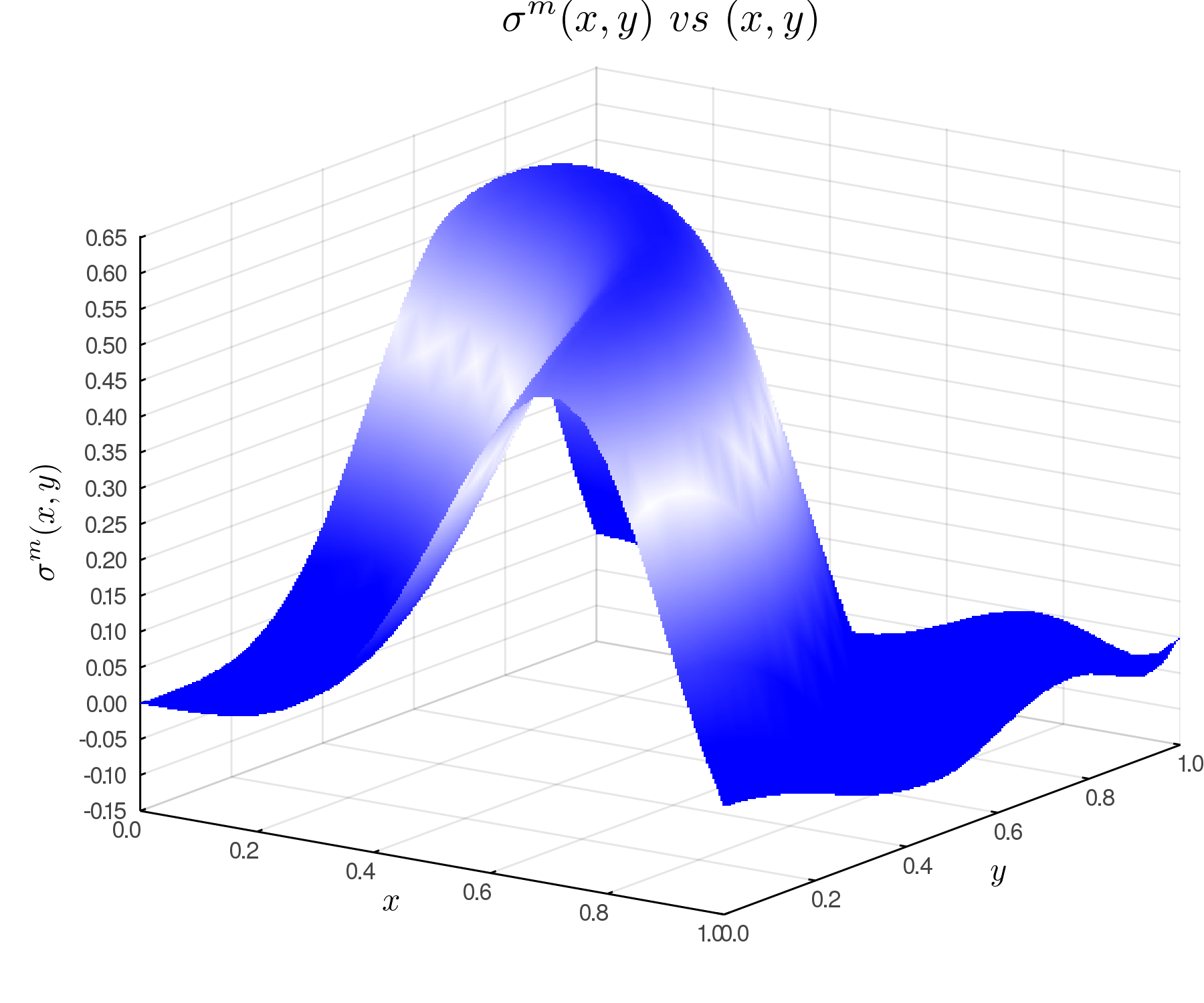}
\caption{\footnotesize Numerical solution for $(\alpha,\beta)=(0.7,1)$.}
\end{subfigure}
\begin{subfigure}[c]{0.33\linewidth}
\includegraphics[width=\linewidth,height=0.84\linewidth]{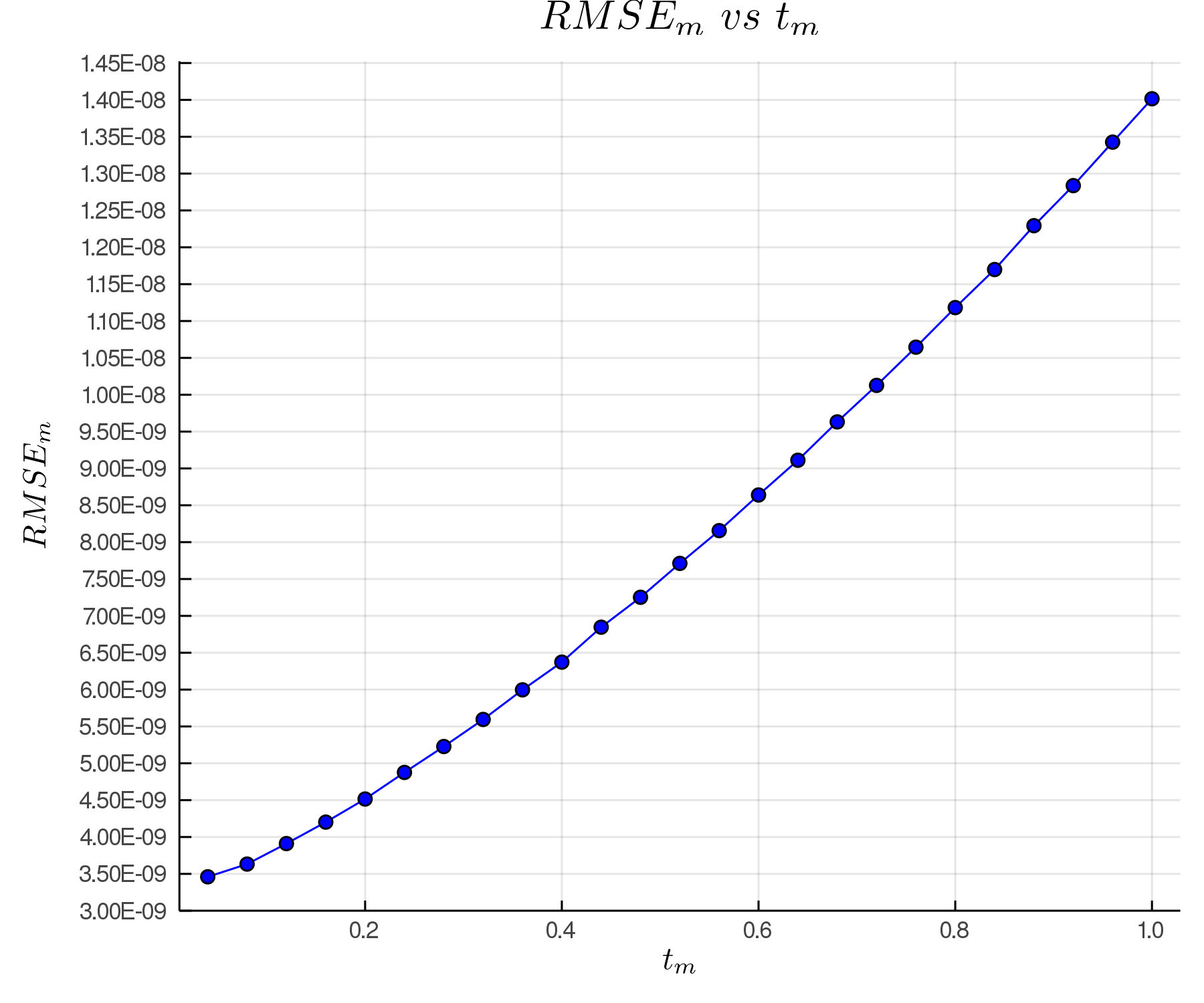}
\caption{\footnotesize $RMSE_m$ for $(\alpha,\beta)=(0.7,1)$.}
\end{subfigure}
\\
\begin{subfigure}[c]{0.33\linewidth}
\includegraphics[width=\linewidth,height=0.84\linewidth]{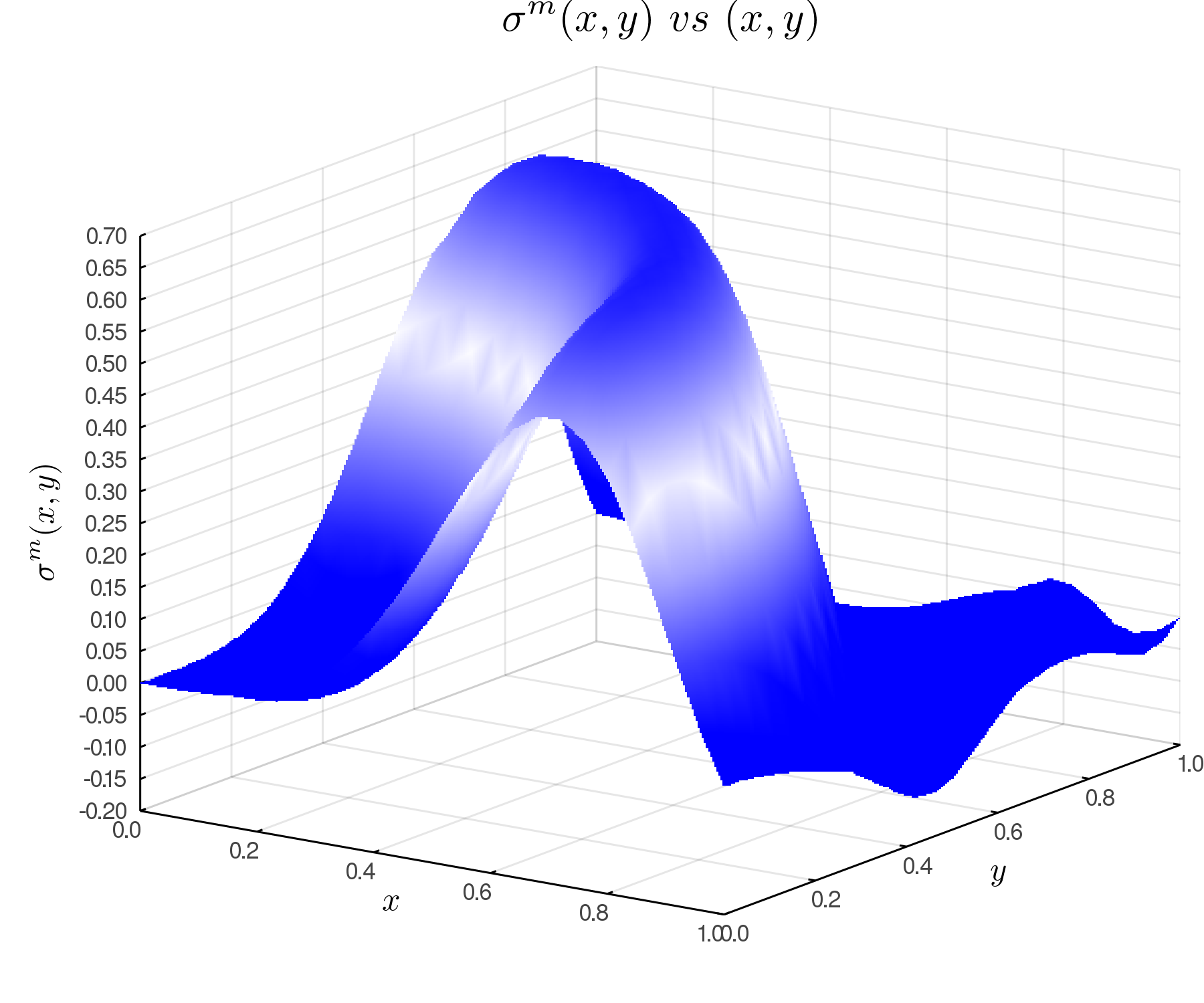}
\caption{\footnotesize Numerical solution for $(\alpha,\beta)=(1,0.75)$.}
\end{subfigure}
\begin{subfigure}[c]{0.33\linewidth}
\includegraphics[width=\linewidth,height=0.84\linewidth]{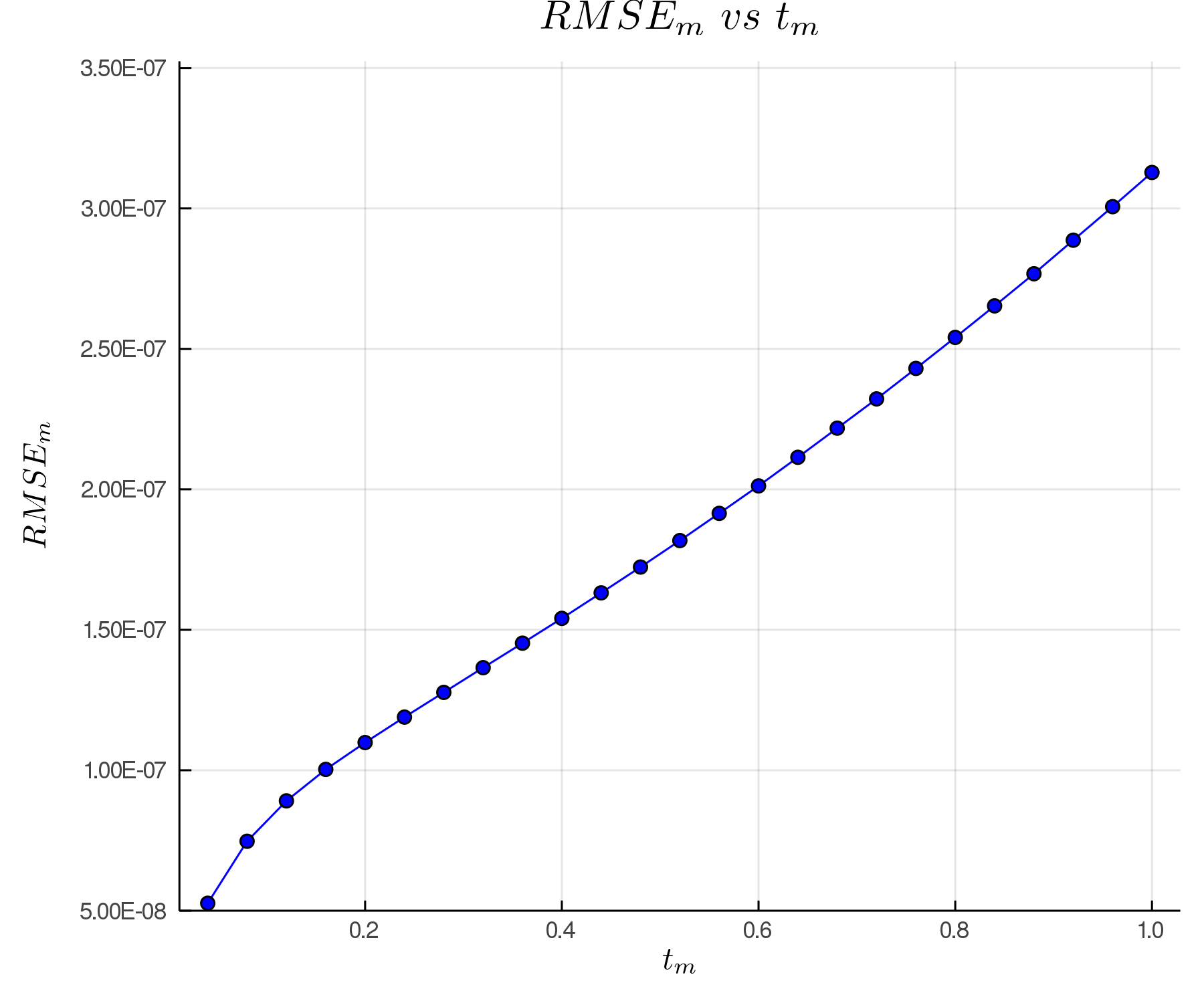}
\caption{\footnotesize $RMSE_m$ for $(\alpha,\beta)=(1,0.75)$.}
\end{subfigure}
\\
\begin{subfigure}[c]{0.33\linewidth}
\includegraphics[width=\linewidth,height=0.84\linewidth]{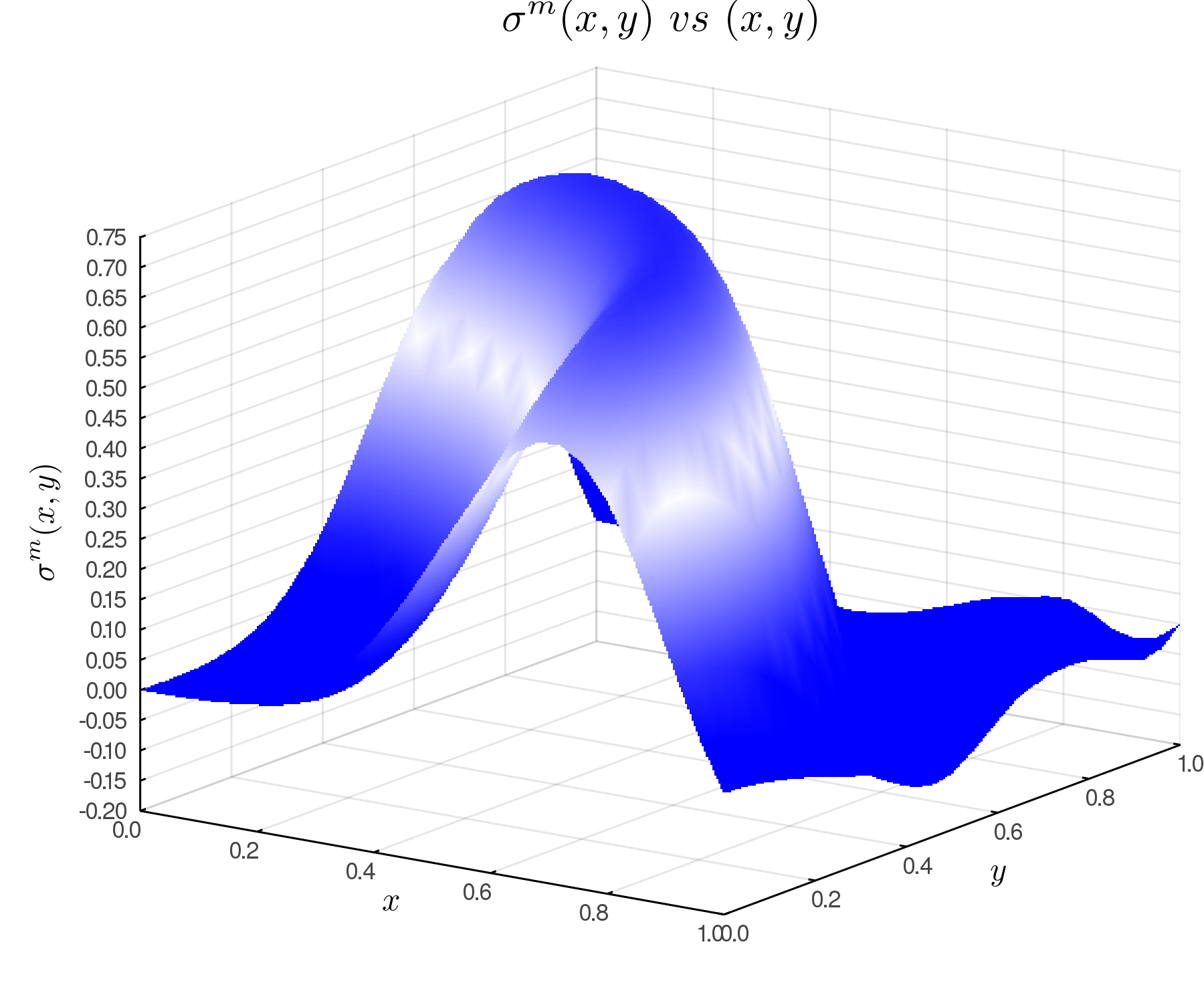}
\caption{ \footnotesize Numerical solution for $(\alpha,\beta)=(0.65,0.8)$.}
\end{subfigure}
\begin{subfigure}[c]{0.33\linewidth}
\includegraphics[width=\linewidth,height=0.84\linewidth]{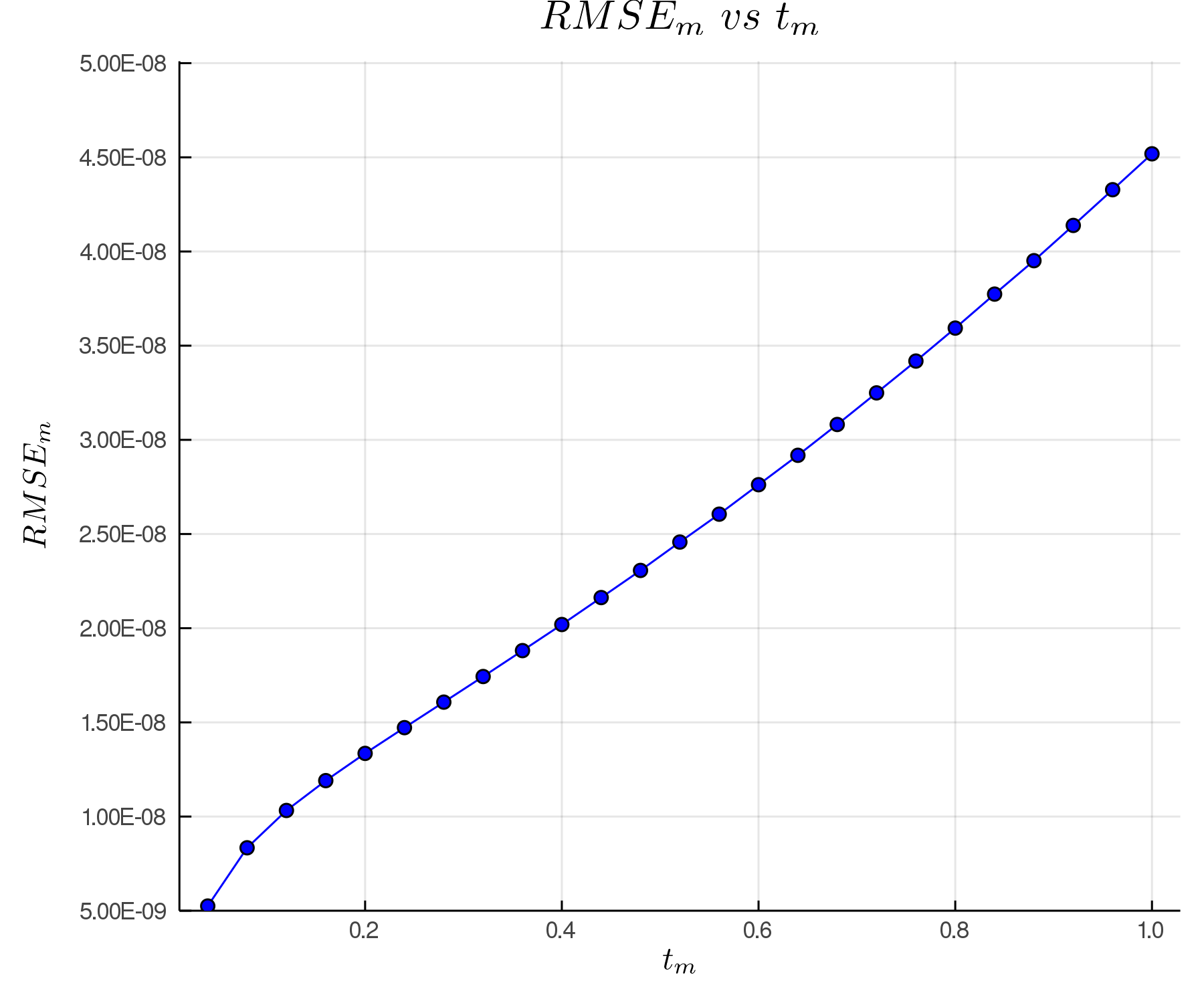}
\caption{\footnotesize $RMSE_m$ for $(\alpha,\beta)=(0.65,0.8)$.}
\end{subfigure}
\caption{The analytical solution and the numerical solutions with respect to space for the final time step are presented. The $RMSE$ is presented with respect to time for the different numerical solutions.}\label{fig:04}
\end{figure}

\begin{figure}[!ht]
\centering
\begin{subfigure}[c]{0.33\linewidth}
\includegraphics[width=\linewidth,height=0.65\linewidth]{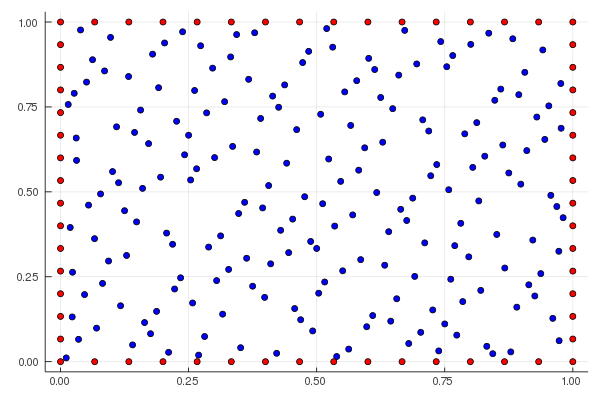}
\caption{$N_p=256$.}
\end{subfigure}
\begin{subfigure}[c]{0.33\linewidth}
\includegraphics[width=\linewidth,height=0.65\linewidth]{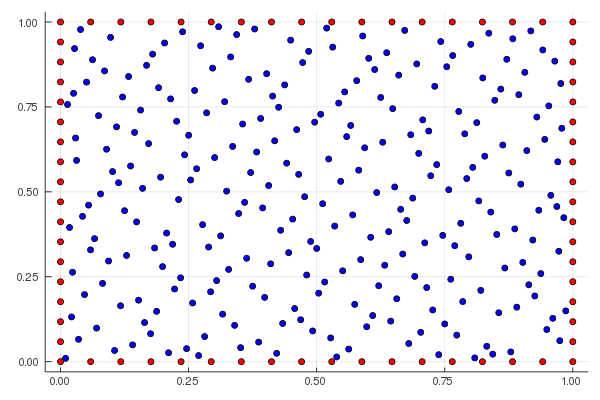}
\caption{$N_p=324$.}
\end{subfigure}
\begin{subfigure}[c]{0.33\linewidth}
\includegraphics[width=\linewidth,height=0.65\linewidth]{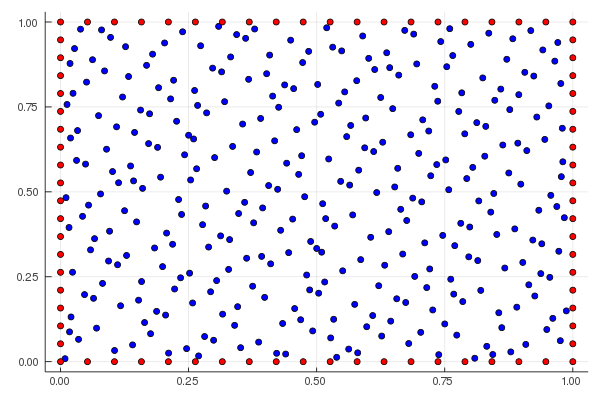}
\caption{$N_p=400$.}
\end{subfigure}
\caption{Different numbers of nodes used.}\label{fig:03}
\end{figure}

\begin{figure}[!ht]
\centering
\begin{subfigure}[c]{0.33\linewidth}
\includegraphics[width=\linewidth,height=0.65\linewidth]{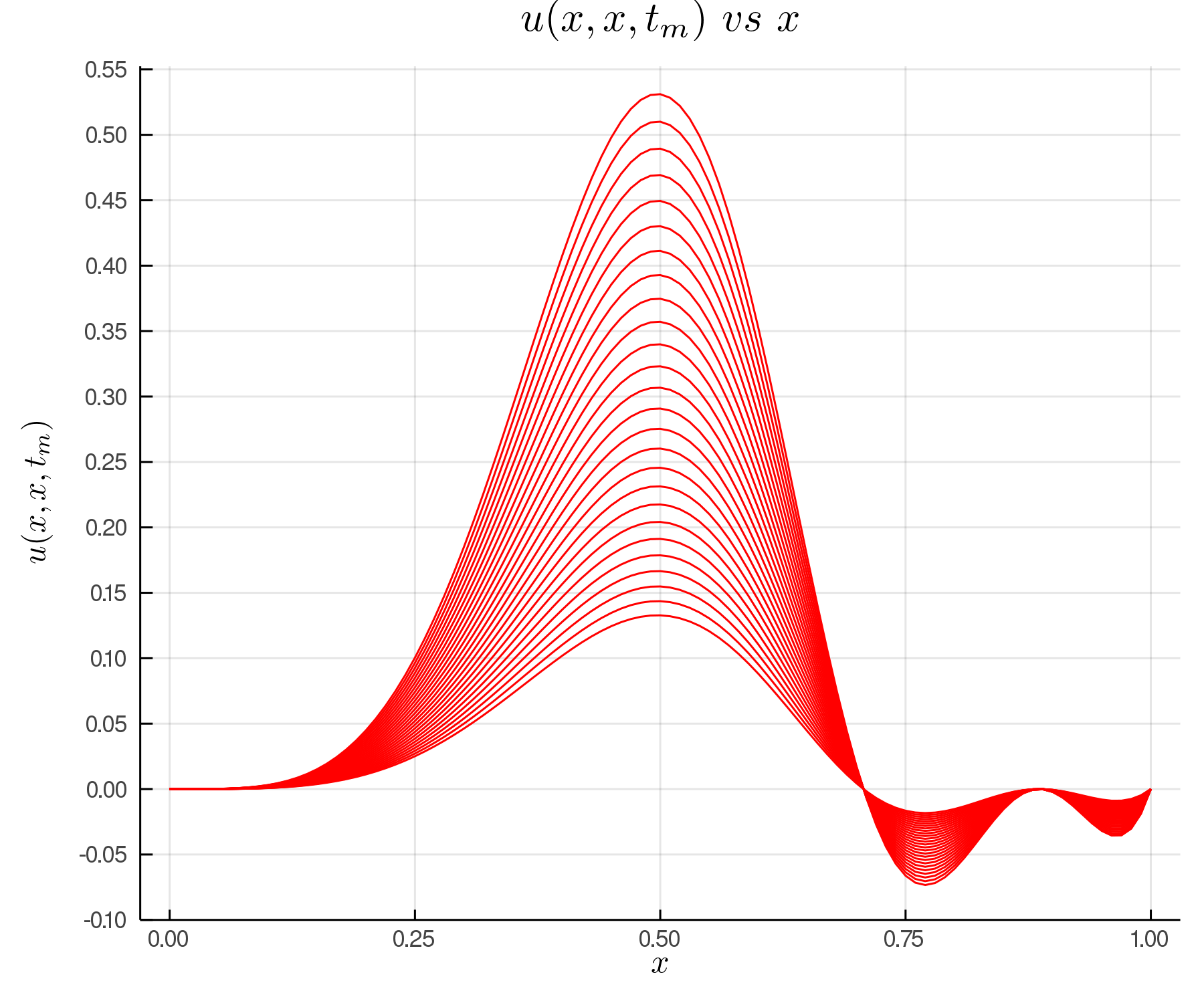}
\caption{\footnotesize Analytical solution.}
\end{subfigure}
\begin{subfigure}[c]{0.33\linewidth}
\includegraphics[width=\linewidth,height=0.65\linewidth]{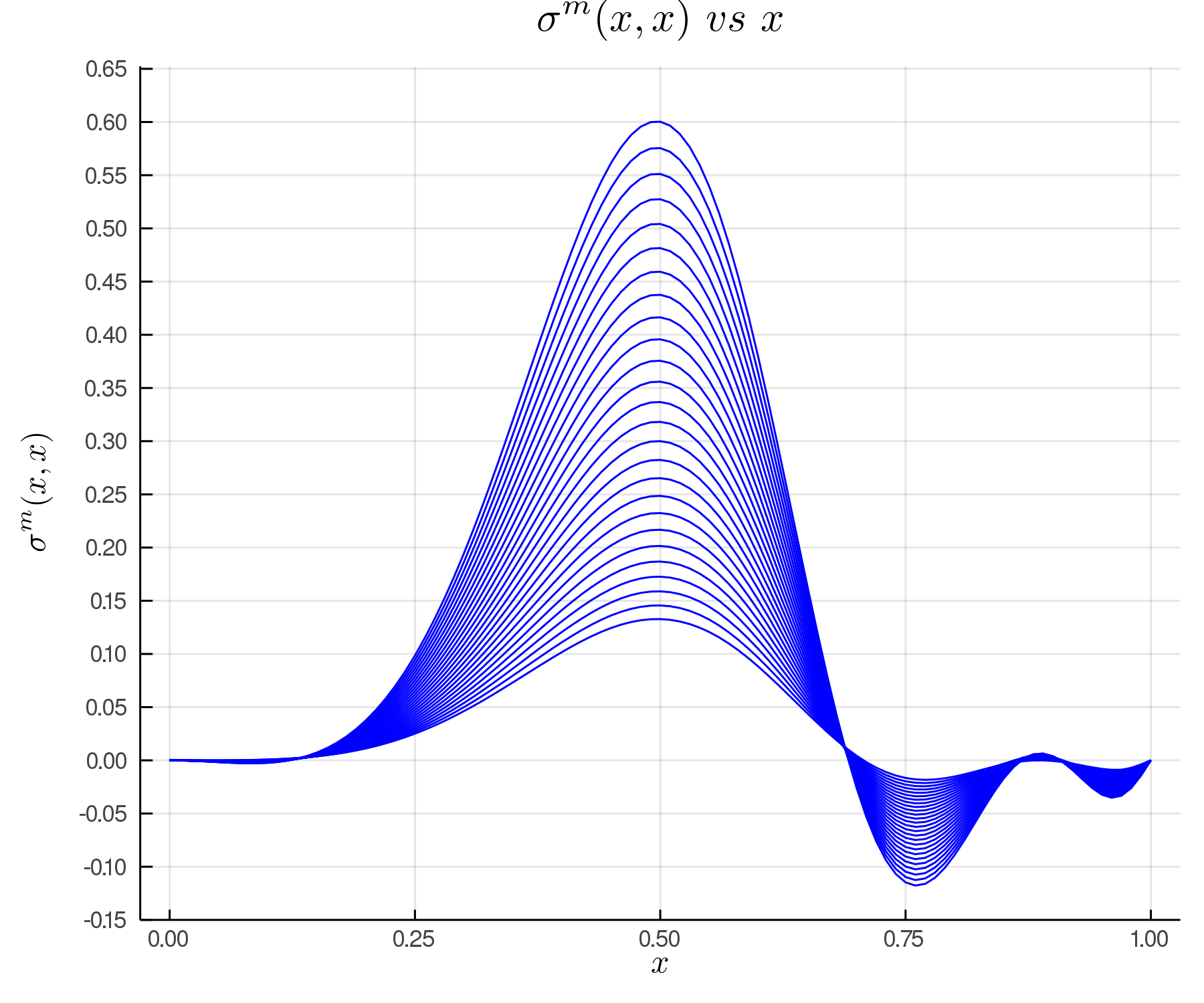}
\caption{\footnotesize Numerical solution for $(\alpha,\beta)=(1,1)$.}
\end{subfigure}
\begin{subfigure}[c]{0.33\linewidth}
\includegraphics[width=\linewidth,height=0.65\linewidth]{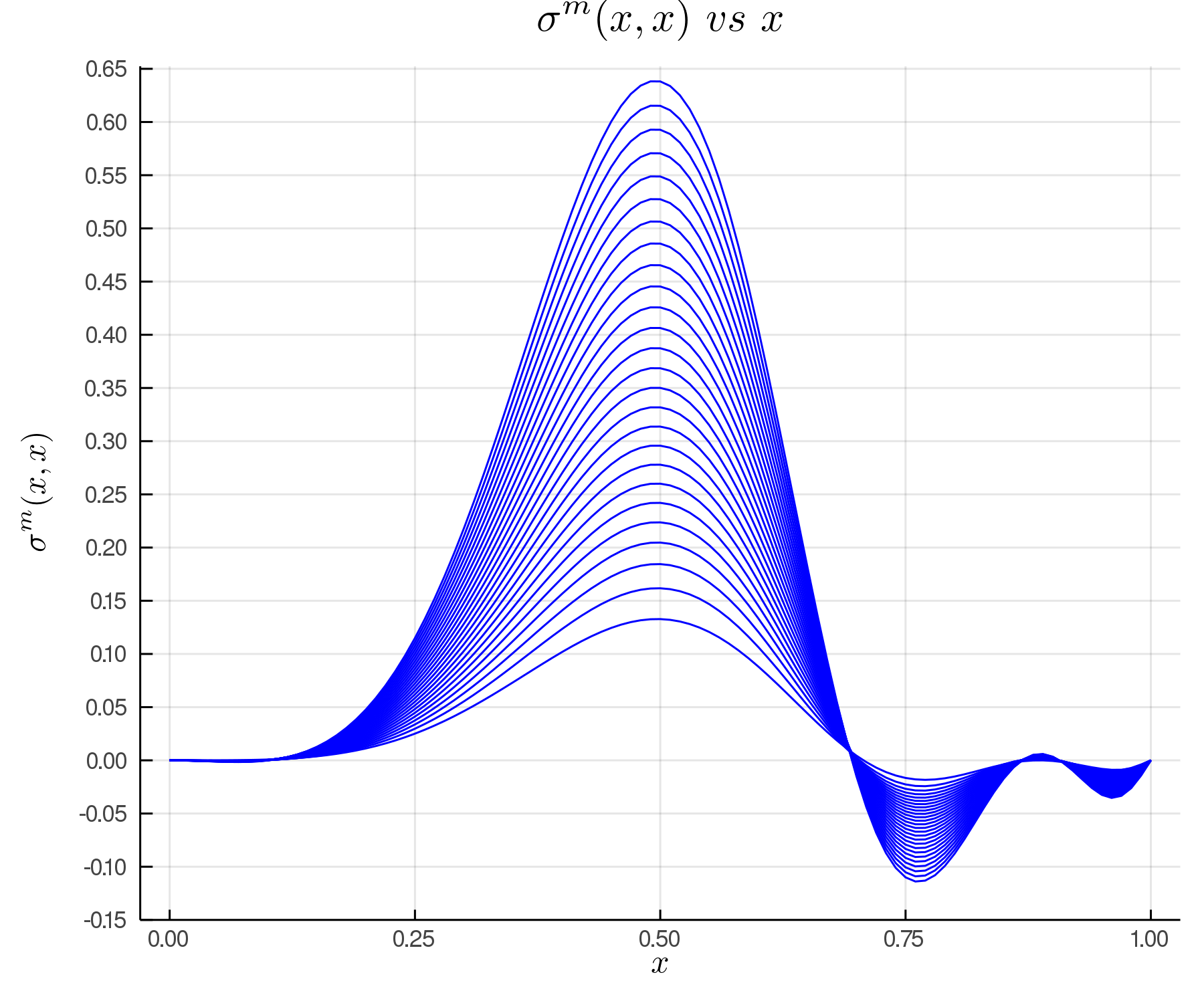}
\caption{\footnotesize Numerical solution for $(\alpha,\beta)=(0.7,1)$.}
\end{subfigure}
\\
\begin{subfigure}[c]{0.33\linewidth}
\includegraphics[width=\linewidth,height=0.65\linewidth]{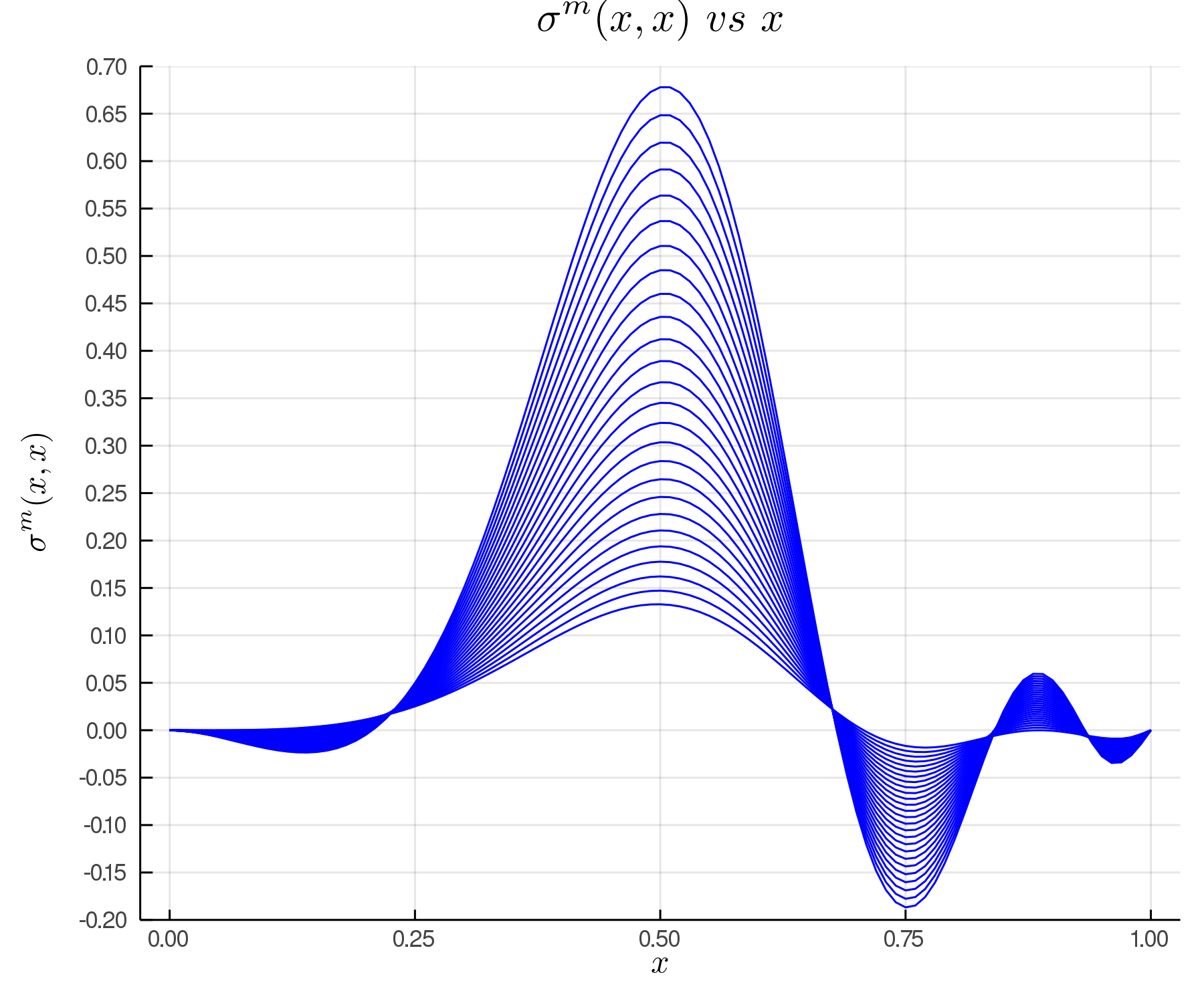}
\caption{\footnotesize Numerical solution for $(\alpha,\beta)=(1,0.75)$.}
\end{subfigure}
\begin{subfigure}[c]{0.33\linewidth}
\includegraphics[width=\linewidth,height=0.65\linewidth]{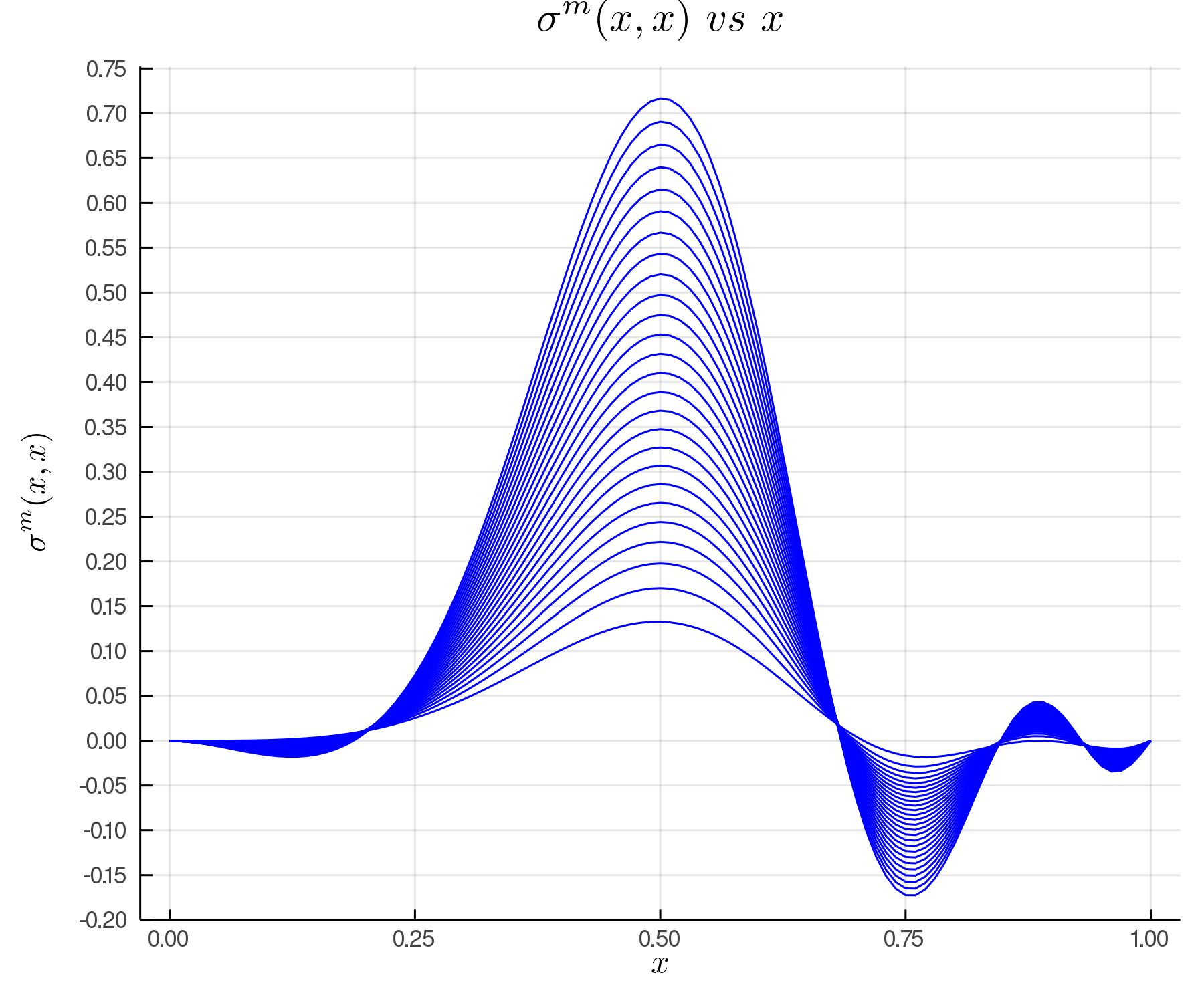}
\caption{ \footnotesize Numerical solution for $(\alpha,\beta)=(0.65,0.8)$.}
\end{subfigure}
\begin{subfigure}[c]{0.33\linewidth}
\end{subfigure}
\caption{The analytical solution and the numerical solutions with respect to space, with $y=x$,  for different moments in time are presented. }\label{fig:05}
\end{figure}

\begin{table}[!ht]
\centering
$
\small
\begin{array}{cccccc}
\toprule
\alpha&\beta&N_p& \cond(G_{\alpha,\beta})&\cond(\widetilde{ G}_{\alpha,\beta})&RMSE\\
\midrule
          &       & 256   & 3.70947E+07 & 2.68764E+00 & 2.99897E-08 \\
    1     & 1     & 324   & 6.33937E+07 & 2.66676E+00 & 1.47561E-08 \\
          &       & 400   & 1.13978E+08 & 2.27838E+00 & 5.94686E-08 \\ \midrule
          &       & 256   & 1.44181E+07 & 1.79504E+00 & 8.35786E-09 \\
    0.7   & 1     & 324   & 2.48632E+07 & 1.71344E+00 & 8.35289E-09 \\
          &       & 400   & 4.54835E+07 & 1.54992E+00 & 1.40151E-08 \\ \midrule
          &       & 256   & 5.21246E+07 & 3.12264E+00 & 1.40965E-07 \\
    1     & 0.75  & 324   & 8.12623E+07 & 3.75596E+00 & 1.97783E-07 \\
          &       & 400   & 1.41019E+08 & 3.16718E+00 & 3.12793E-07 \\ \midrule
          &       & 256   & 1.42528E+07 & 2.13174E+00 & 9.01938E-09 \\
    0.65  & 0.8   & 324   & 2.42863E+07 & 2.06023E+00 & 1.28433E-08 \\
          &       & 400   & 4.37752E+07 & 1.80894E+00 & 4.51864E-08 \\ \bottomrule
\end{array}
$
\caption{Values obtained for the different numerical solutions, the value of $RMSE$ is presented for the final time step.}\label{tab:02}
\end{table}

\end{example}

The errors in Figures \ref{fig:02} and \ref{fig:04} show an increasing behavior with time, which is consistent with the  condition \eqref{eq:S3-02}. For the case where $\alpha=1$, the errors fulfill the following condition

\begin{eqnarray}
\mathcal{O}_{\alpha}^m(x)=\mathcal{O}_{\alpha}^m\left(x,\mathcal{O}_{\alpha}^{m-1}(x) \right),
\end{eqnarray}

however, the condition \eqref{eq:S3-02} is still fulfilling implicitly.  The results obtained in the previous examples could be improved by implementing one or more of the following strategies:

\begin{itemize}
\item[i)] Selecting a smaller $dt$ time step.
\item[ii)] Working with a greater number $N_p$ of nodes.
\item[iii)] Changing the set of radial functions  $\set{\Phi(x,x_j)}_{j=1}^{N_p}$.
\end{itemize}

To keep errors under control, strategy $iii)$ would be the most recommended. Polyharmonic radial functions \cite{wendland} could be used

\begin{eqnarray*}
\Phi(x,x_j)=\norm{x-x_j}_2^{2n+1}, & n\in \nset{N},
\end{eqnarray*}

or multiquadratic radial functions \cite{wendland}

\begin{eqnarray*}
\Phi_{\epsilon}(x,x_j)=\left[1+\left( \epsilon \norm{x-x_j}_2 \right)^2 \right]^{\mu/2}, & \mu \in [-1,1] \setminus \set{0},
\end{eqnarray*}

these last functions incorporate a parameter $\epsilon\in \nset{R}_{>0}$, known as a shape parameter, which being varied allows to improve the errors of the numerical solutions without the need to decrease the time step or increase the number of nodes. However, finding the optimal shape parameter $\epsilon$ for each problem is computationally expensive.

In general, given the expression \eqref{eq:S3-01}, which is a consequence of the memory phenomenon in the fractional differential operator in time, a prudent strategy would be to leave as a last resort, to improve errors in numerical solutions, use radial basis functions with a shape parameter. The latter with the aim of not increase to a large degree the computational cost to solve multidimensional fractional partial differential equation systems.

\section{Conclusions}

In this work, the flexibility of the radial basis functions scheme was shown to solve multidimensional problems with various types of nodes and it was also shown how to reduce the condition number of the matrices involved. Problems related to the space-time-fractional Black-Scholes equations were solved in one and two dimensions, reducing the condition number of the discretization matrices of the differential operator by approximately less than one percent of their original value. Chebyshev nodes were used and also Halton nodes combined with Cartesian nodes, but in general, any distribution of nodes, uniform or non-uniform, and combinations of them can be used.

The easy implementation of the radial basis function scheme to solve fractional equations allows considering different types of placement nodes and generalizing to large dimensions. This allows us to focus on making the scheme more stable and efficient by reducing the condition number of the matrices involved in the process. As shown, the meshless method via radial basis functions is implemented to solve time-space-fractional equations of type Black-Scholes. The results show that, although errors grow over time, is an efficient technique and may be considered as a numerical technique for solving different one-dimensional or multidimensional fractional  partial differential equations

The schemes that use radial basis functions are easy to implement compared to finite element schemes or finite difference schemes, this characteristic becomes more evident when attacking problems in multiple dimensions, as a consequence of the dimensional invariance of the radial basis functions methodology.  However, even with this advantage over finite differences or finite element, before using radial basis function schemes, the computational cost and susceptibility to numerical errors must be considered, since the matrices involved can be analytically invertible but numerically singular.

\bibliography{Biblio}
\bibliographystyle{unsrt}

\nocite{kumar2014numerical}
\nocite{odibat2009computational}
\nocite{liang2010solution}
\nocite{mohammadi2016fractional}
\nocite{yang2010numerical}
\nocite{chen2014recent}
\nocite{uddin2011rbfs}
\nocite{odibat2009variational}

\end{document}